\documentclass{article}

\usepackage[colorlinks	= true,
			raiselinks	= true,
			linkcolor	= MidnightBlue,
			citecolor	= Mahogany,
			urlcolor	= ForestGreen,
			pdfauthor	= {Peyman Mohajerin Esfahani},
			pdftitle	= {},
			pdfkeywords	= {},   
			pdfsubject	= {},
			plainpages	= false]{hyperref}
\usepackage[usenames, dvipsnames]{xcolor}
\usepackage{dsfont,amssymb,amsmath,graphicx,fancyhdr,mdframed}
\usepackage{amsfonts,dsfont,mathtools, mathrsfs,amsthm,wrapfig} 
\usepackage{mdframed,ragged2e,subfigure}
 
\usepackage{enumitem}
\usepackage{constants}
\allowdisplaybreaks
\usepackage{algorithm}
\usepackage[noend]{algpseudocode}

\allowdisplaybreaks
\date{\today}

\newtheorem{thm}{Theorem}[section]
\newtheorem{prop}[thm]{Proposition}
\newtheorem{defn}[thm]{Definition}

\newtheorem{lem}[thm]{Lemma}

\newtheorem{assumption}[thm]{Assumption}

\newtheorem{remark}[thm]{Remark}
\newtheorem{corr}[thm]{Corollary}

\def \l{\left}
\def \r{\right}
\def \a{\alpha}
\def \R{\mathbb{R}}
\def \Rfrak{\mathfrak{r}}
\def \wh{\widehat}
\def \E{\mathbb{E}}

\def \ol{\overline}
\def \b{\beta}
\def \B{\mathrm{B}}
\def \g{\gamma}

\def \tr{\operatorname{trace}}
\def \vol{\text{d}}

\def \W{\mathrm{W}}

\def \der{\operatorname{d}\hspace{-1pt}}

\def \P{\mathcal{P}}

\def \P{\mathfrak{p}}

\def \dist{\operatorname{dist}}

\DeclareMathOperator*{\argmin}{arg\,min}

\def \b{\beta}

\def\Q{\mathfrak{q}}

\def\CC{\mathrm{K}}

\def\MY{\textup{MY}}

\def\FB{\textup{FB}}
\def\vol{\textup{vol}}
\def\K{{K}}

\newconstantfamily{fcn}{symbol=c}

\title{The Forward-Backward Envelope for Sampling \\
with the Overdamped Langevin Algorithm}

\author{Armin Eftekhari \and Luis Vargas\and Konstantinos Zygalakis\thanks{The authors are ordered alphabetically. AE is with the Department of Mathematics and Mathematical Statistics, Umea University. LV and KZ are with the the School of Mathematics, University of Edinburgh, and the Maxwell Institute for Mathematical Sciences.}}

\begin{document}

\maketitle

\begin{abstract}
    In this paper, we analyse a proximal method based on the idea of forward-backward splitting for sampling  from distributions {with} densities {that} are not necessarily smooth. In particular, we study the non-asymptotic properties of the Euler-Maruyama discretization of the Langevin equation, where the forward-backward envelope is used to deal with the non-smooth part of the dynamics. An  advantage of this envelope, when compared to widely-used  Moreu-Yoshida one and the MYULA algorithm, is that it maintains the MAP estimator of the original non-smooth distribution. We 
    {also}
    study a number of 
    numerical experiments  
    {that support}
    our theoretical findings.
\end{abstract}

\section{Introduction}


The problem of calculating expectations with respect to a 
probability distribution {$\P$ in $\R^d$} is ubiquitous  throughout applied mathematics, 
statistics, molecular dynamics, statistical physics and other 
fields. {In practice, often}
{$d$ is large, which} renders 
deterministic techniques, such as quadrature {methods}, {computationally intractable.}  {In contrast,} probabilistic methods do not suffer {from} the 
curse of dimensionality and are often the method of choice
{when the dimension $d$ is large.}
{In particular,}
Markov {c}hain Monte Carlo 
(MCMC) {methods} 
{are} based on the construction of a Markov chain 
{in} $\mathbb{R}^{m}$ 
{with} $m \geq d$, for which the 
{invariant} distribution 
(or 
{its} suitable marginal)
coincides with 
{the target distribution}~$\P$~\cite{brooks2011handbook}. 

{Often, such} Markov chains are based on the  discretization of stochastic 
differential equations (SDEs).  One
such SDE, which 
is also the 
focus {of} this paper, is the
(overdamped) Langevin equation 
\begin{equation} \label{eq:Langevin}
\der X_{t}=-\nabla f(X_{t})\der t+\sqrt{2}\der W_{t},   
\end{equation}
where $\{W_t\}_t$ is the standard $d$-dimensional Brownian motion and $\nabla f$ denotes the gradient of a continuously-differentiable function $f:\R^d\rightarrow\R$.
Under mild assumptions on $f$, 
the 
dynamics of~\eqref{eq:Langevin} are ergodic with respect to the distribution
$\P \propto e^{-f}$.
{In particular, $\P$ is the invariant distribution of~\eqref{eq:Langevin}~\cite{MT07}.}

{The discretization of~\eqref{eq:Langevin}, however, requires special care,}
since 
the {resulting discrete Markov chain}
might 
not be ergodic \cite{RT96}. 
In addition, even if 
ergodic, the {resulting discrete Markov} chain {often has a different invariant distribution than~$\P$, known as the numerical invariant distribution~$\widehat{\P}$.} The study of the asymptotic error between the 
numerical {invariant} distribution~$\widehat{\P}$ and
{the target distribution~$\P$}
has 
received 
{considerable} attention {recently}~\cite{MST10,AGZ14}. {In particular,}~\cite{MST10} 
{investigated}
the effect of discretization on the
convergence of the ergodic averages, and \cite{AGZ14} 
{presented} general order conditions 
{to ensure that}
the numerical invariant 
{distribution}
{accurately} approximates the
{target distribution.}

Another {active} line of research 
{quantifies}
the nonasymptotic error between the 
numerical {invariant distribution}~$\widehat{\P}$
and the 
{target distribution $\P$.}
{In particular,} 
{when $\P$ is a}
{smooth}
and strongly log-concave 
 distribution,~\cite{D17b} 
established non-asymptotic  bounds in total variation distance 
for the Euler-Maruyama {discretization of \eqref{eq:Langevin}, commonly known as the unadjusted Langevin algorithm (ULA).}
These results 
have {also} been extended to the Wasserstein distance $\W_{2}$ in 
\cite{DM17,D17a,DM19,DK19,DMM18}, to name a few.
{Typically, these works study} the number of iterations 
that the {numerical}
integrator 
{would} require to achieve a 
{desired} accuracy, 
when applied to {a} 
target {distribution $\P$} with {a known} condition number. 

{In fact, the above strong log-concavity of $\P$  can be substantially relaxed. In particular, using a variant of the reflection coupling, the recent work~\cite{eberle2019quantitative} derived non-asymptotic bounds for the ULA in the Wasserstein distance~$\W_1$, when $\P$ is strictly log-concave \emph{outside} of a ball in $\R^d$. Similar results for the Wasserstein distance~$\W_2$ have also been presented in~\cite{majka2020nonasymptotic}. }

Within the class of log-concave distributions, a significant 
challenge for the Langevin {diffusion in~\eqref{eq:Langevin}}
arises when the target 
{distribution}
$\P$ is not smooth and/or {has a compact (convex) support in $\R^d$.} 
One {approach to address}
{this challenge is to replace the non-smooth 
distribution~$\P$ with a smooth proxy obtained via  the so-called Moreu-Yoshida (MY) envelope. This new smooth density remains log-concave and, hence, {amenable to the non-asymptotic results discussed earlier.} When the support of $\P$ is also compact,  proximal Monte Carlo methods  have been explored in~\cite{brosse2017sampling,P16,DMP18}. It is also worth noting that~\cite{bubeck2018sampling} pursued a different approach for sampling from compactly-supported densities that does not involve the MY envelope.}

A potential drawback of the above approach is that the MY envelope  
often
 does not maintain the maximum a posteriori (MAP) estimator.
 That is, the above approach
alters the 
location at which the new (smooth) density 
{reaches}
its maximum. This {is a well-known issue}
in the {context of }(non-smooth) convex optimization 
and {is often resolved by appealing to}
the proximal gradient method. 
The latter can be understood as the Euler discretization of the so-called forward-backward (FB) envelope~{\cite{stella2017forward}}. 

\paragraph{Contributions.} 
This work explores and analyzes the 
use of the FB envelope for sampling from 
non-smooth and compactly-supported log-concave 
distributions. 
In analogy with the Langevin proximal Monte Carlo, we replace the 
non-smooth density 
{with}
a smooth {proxy obtained via the FB envelope.} In particular, this proxy is strongly log-concave over long distances. 

Crucially, the new proxy also
maintains the 
MAP estimator {, under certain assumptions}. {However,}  this {improvement} comes at the cost of 
requiring additional smoothness for the smooth 
 part of the density.
Lastly, the strong convexity of the new proxy over long distances allows us to utilise the work of \cite{eberle2019quantitative} 
{to obtain}
non-asymptotic guarantees for our method in the {Wasserstein distance $\W_1$}.

{In addition to investigating the use of FB envelope in sampling, this work has the following contributions:
\begin{itemize}
    \item It introduces a general theoretical framework for sampling from non-smooth densities by introducing the notion of admissible envelopes. MY and FB envelopes are both instances of admissible envelopes. 
    \item It proposes a new Langevin algorithm to sample from non-smooth densities, dubbed EULA, which generalizes MYULA. EULA can work with any admissible envelope (e.g., MY or FB) and can handle a family of increasingly more accurate envelopes rather than a fixed envelope.   
\end{itemize}
}

\paragraph{{Organization.}}
The rest of the paper is organised as follows. Section~\ref{sec:statement} formalizes the problem of sampling from a non-smooth and compactly-supported log-concave distribution. As a proxy for this non-smooth distribution, its (smooth) Moreau-Yosida (MY) envelope is reviewed in Section~\ref{sec:my-Env}. This section also explains the main limitation of MY envelope, i.e., its inaccurate MAP estimation. 
In Section~\ref{sec:FB}, we introduce the forward-backward (FB) envelope which overcomes the key shortcoming of the MY envelope.

Section~\ref{sec:newALg} introduces and analyses EULA, an extension of the popular ULA for sampling from a non-smooth distribution. EULA  can be adapted to various envelopes. In particular, MYULA from~\cite{brosse2017sampling} is a special case of EULA for the MY envelope. Section~\ref{sec:proofMain} proves the iteration complexity of EULA and Section~\ref{sec:numerics} presents a few numerical examples to support the theory developed here.

\section{Statement of the Problem}\label{sec:statement}

Consider a compact convex set $\CC\subset\R^d$. For a pair of functions $\ol{f}:\R^d\rightarrow\R$ and $\ol{g}:\R^d\rightarrow\R$, our objective in this work is to sample from the probability distribution 
\begin{equation}
\P(x) :=  
\begin{cases}
\frac{e^{-\ol{f}(x)-\ol{g}(x)}}{\int_{\CC} e^{-\ol{f}(z)-\ol{g}(z)} \,\der z } & x\in \CC \\
0 & x\notin \CC,
\end{cases}
\label{eq:pi0}
\end{equation}
whenever the ratio above is well-defined. In order to sample from $\P$, we only have access to the gradient of $\ol{f}$ and the proximal operator for $\ol{g}$, to be defined later. Our assumptions on $\CC,\ol{f},\ol{g}$ are detailed below.


\begin{assumption}\label{assumption:fg}We make the following assumptions:

\begin{enumerate}[leftmargin=*,label={(\roman*)}]
\item\label{item:CAssumption} For radii $R\ge r>0$, assume that  $\CC\subset \R^d$ is a  compact convex body that satisfies $\B(0,r)\subset \CC \subset \B(0,R)$. Here, $\B(0,r)$ is the Euclidean ball of radius $r$ centered at the origin. 
\item\label{item:fAssump} Assume also that $\ol{f}:\R^d\rightarrow\R $ is a convex function that is three-times continuously differentiable. 
\item Assume lastly that  $\ol{g}:\R^d\rightarrow(-\infty,\infty]$ is a proper closed convex function.
Moreover, we assume that~$\ol{g}$ is continuous.\footnote{In Assumption~\ref{assumption:fg}\ref{item:fAssump}, the requirement that $\ol{g}$ is a proper closed convex function  implies that $\ol{g}$ is lower semi-continuous, but not necessarily continuous~\cite{bertsekas2009convex}. The latter stronger requirement of continuity for $\ol{g}$ is needed in this work.}
\end{enumerate}
\end{assumption}

A few important remarks about Assumption \ref{assumption:fg} are in order.  
{First, in the special case when $\ol{f}$ is a convex quadratic~\cite{LFC21}, the assumption of thrice-differentiability above is trivially met and some of the developments below are simplified. However, our more general setup here necessitates the thrice-differentiability above and results below in more involved technical derivations.}

Second, instead of the two functions $\ol{f},\ol{g}$, it will be more convenient to work with two new functions~$f,g$, without any loss of generality. 
More specifically, consider a convex function $f$ that coincides with $\ol{f}$ on the set $\CC$, has a compact support and a continuously differentiable Hessian.

For this function~$f$, the compactness of $\CC$ and smoothness of $\ol{f}$  in Assumption~\ref{assumption:fg} together imply that $f,\nabla f,\nabla^2 f$ are all Lipschitz-continuous functions. To summarize, for the function~$f$ described above, there exist nonnegative constants $\lambda_0,\lambda_1,\lambda_2,\lambda_3$ such that 
 \begin{subequations} 
\begin{align}
& f(x) = 0, \qquad \text{if } \|x\|_2 \ge \lambda_0,\label{eq:compactSuppf} \\
& |f(x)-f(y)| \le \lambda_1 \|x-y\|_2, \qquad x,y\in \R^d, \label{eq:fLipsch} \\
& \|\nabla f(x) - \nabla f(y) \|_2 \le \lambda_2 \|x-y\|_2, \qquad x,y\in \R^d,\label{eq:gradLipschitz} \\
      & \|\nabla^2 f(x) - \nabla^2 f(y) \| \le \lambda_3 \|x-y\|_2,
     \qquad  x,y\in \R^d. 
     \label{eq:hessianLipschitz}
\end{align}
\label{eq:lipsAll}\end{subequations}Let us also define the proper closed convex function 
\begin{equation}
g := \ol{g}+1_\CC,
\label{eq:g_gbar}
\end{equation}
where $1_\CC$ is the indicator function for the set $\CC$. That is, $1_\CC(x)=0$ if $x\in \CC$ and $1_\CC(x)=\infty$ if $x\notin \CC$. 
The compactness of $\CC$ and continuity of~${\bar{g}}$ in Assumption~\ref{assumption:fg} together imply  that $g$ is finite, when its domain is limited to the set $\CC$. Outside of the set~$\CC$, $g$ is infinite. To summarize, the new function $g$ is lower semi-continuous and also satisfies 
\begin{align}
    \max_{x\in \CC}|g(x)| <\infty, \qquad
     g(x)=\infty, \quad \text{if }x\notin \CC. 
    \label{eq:propsG}
\end{align}
We can now revisit  \eqref{eq:pi0} and, using the new functions $f,g$, we rewrite the definition of $\P$ as 
\begin{subequations}\begin{align}
& \P(x) :=  \frac{e^{-F(x)}}{\int_{\R^d} e^{-F(z)} \,\der z }, \qquad x\in \R^d, \label{eq:pi}\\
& F(x):=\ol{f}(x)+\ol{g}(x)+1_\CC(x)= f(x) + g(x),
\qquad x\in \R^d.
\label{eq:BigF}
\end{align}\end{subequations}Above, $F$ is often referred to as the potential associated with $\P$. The last identity above holds by construction. Indeed, on the set $\CC$, the functions $f$ and $\ol{f}$ coincide. Likewise, on the set $\CC$,  the functions $g$ and $\ol{g}$ coincide. On the other hand, outside of the set $\CC$, both sides of the last equality above are infinite. 

In view of \eqref{eq:BigF}, we will often use $f$ and $\ol{f}$ interchangeably throughout this work, depending on the context. Likewise, we will use $g$ and $\ol{g}$ interchangeably.  Note also that the integral in the denominator above is finite by Assumption \ref{assumption:fg}. When there is no confusion, we will overload our notation and use $\P$ to also denote the probability measure associated with the law $\P$.

Since $g$ is not differentiable, $F$ in \eqref{eq:BigF} is itself non-differentiable. In turn, this means that one cannot use gradient based algorithms such as ULA to sample from  $\P\propto e^{-F}$~\cite{DK19}. One way to deal with this issue is to replace $F$ with a smooth function $F_{\gamma}$, which we will refer to as an envelope, to which we can  {then} apply ULA.  It is reasonable to  require this envelope $F_\gamma$ to fulfill the following admissibility assumptions.


\begin{defn}[{\sc Admissible envelopes}]\label{assump:envAssumption}
For $\g^0>0$, the functions $\{F_\gamma:\R^d\rightarrow[-\infty,\infty]: \g\in (0,\g^0)\}$ are admissible envelopes of $F$ if
\begin{enumerate}[leftmargin=*,label={(\roman*)}]
\item\label{item:integrable} There exists a  function $F^0:\R^d\rightarrow[-\infty,\infty]$ such that $e^{-F^0}$ is integrable, and $F_\g$ dominates $F^0$. That is, 
    \begin{equation*}
        \int_{\R^d} e^{-F^0(z)} \, \der z<\infty, \qquad F_\g(x) \ge F^0(x), \qquad x\in \R^d, \, \g \in(0, \g^0). 
    \end{equation*}
      \item\label{item:pntCvg} $F_\gamma$  converges pointwise to $F$, i.e., 
    $
        \lim_{\gamma\rightarrow0}F_\gamma(x)= F(x)$ for every $x\in \R^d$.
    
    \item\label{item:gradLips} $F_\g$ is $\lambda_\g$-smooth, i.e., there exists a constant $\lambda_{\gamma}\ge 0$ such that 
     \begin{equation*}
         \|\nabla F_\gamma(x)-\nabla F_\gamma(y)\|_2\le \lambda_{\gamma} \|x-y\|_2, \qquad x,y\in \R^d, \, \g\in (0,\g^0). 
     \end{equation*}

\end{enumerate}
\end{defn}
If $\{F_\g:\g \in (0, \g^0)\}$ are  admissible envelopes of $F$, we can define the corresponding probability densities
\begin{equation}
    \P_\gamma(x) := \frac{e^{-F_\gamma(x)}}{\int_{\R^d}e^{-F_\gamma(z)}\, \der z}, \qquad x\in \R^d, \, \g \in (0,\g^0). 
    \label{eq:piGamma}
\end{equation}

\begin{remark}\label{rem:imp}
In the definition above, the property \ref{item:integrable} implies that $\P_\gamma$  can be normalized. This observation, combined with the property \ref{item:pntCvg}, imply {after an application of the dominated convergence theorem} that 
\begin{equation*} 
    \lim_{\g\rightarrow 0}\P_\g(x)  = \P(x), \qquad x\in \R^d. 
\end{equation*}
where  the probability measures $\P$ and $\P_\gamma$ are  defined  in~\eqref{eq:pi} and \eqref{eq:piGamma}, respectively. That is, $\P_\g$ converges weakly to $\P$ in the limit of~$\g\rightarrow 0$. (For completeness, the proof of this claim is included in the appendix.)  
In other words, we can use $\P_\g$ as a proxy for $\P$, provided that $\g$ is sufficiently small. 
Finally, as we will see shortly, the property~\ref{item:gradLips} guarantees the convergence of the ULA to an  invariant distribution close to~$\P_\g$, provided that the  step size of the ULA is small. 
\end{remark}




\section{Moreau-Yosida Envelope and Its Limitation \label{sec:my-Env}}

For $\g>0$, let us define 
\begin{equation}
    F^{\MY}_\g(x) := f(x) + g_{\g}(x), \qquad x\in \R^d,
    \label{eq:MYEnv}
    \tag{MY}
\end{equation}
where 
\begin{equation}\label{eq:MY}
g_{\gamma}(x):=\min_{z \in \mathbb{R}^{d}} \left\{g(z)+\frac{1}{2\gamma}\|x-z\|_2^{2} \right\}
\end{equation}
is the Moreau-Yosida (MY) envelope of $g$. Somewhat inaccurately, we will also refer to $F_\g^\MY$ as the MY  envelope of $F$, to distinguish $F^\MY_\g$ from its newer alternatives.
It is well-known that $g_\g$ is $\g^{-1}$-smooth and that $g_\g$ converges pointwise to $g$ in the limit of $\g\rightarrow0$.   
These facts enable us to establish the admissibility of MY envelopes, as detailed below. {All proofs are deferred to the appendices. We note that the result below closely relates to \cite[Proposition 1]{DMP18}.}
\begin{prop}[{\sc Admissibility of MY envelopes}]\label{prop:my} Suppose that Assumption \ref{assumption:fg} is fulfilled. Then $\{F_\g^{\MY}:\g> 0\}$ are admissible envelopes of $F$ in \eqref{eq:BigF}. In particular, $\nabla F_\g^{\MY}$  is $(\lambda_2+\g^{-1})$-Lipschitz continuous, and given by the expression
\begin{align}
    & \nabla F_\g^{\MY}(x) = \nabla f(x) + \frac{x- P_{\g g}(x)}{\g} \in \nabla f(x) + \partial g(P_{\g g}(x)),  \nonumber \\
    & P_{\g g}(x) := \argmin_{z\in \R^d} \l\{ g(z) + \frac{1}{2\g}\|x-z\|_2^2\r\}.
    \label{eq:gradMY+proximal}
\end{align}
Above, $\lambda_2$ was defined in \eqref{eq:gradLipschitz},  $P_{\g g}:\R^d\rightarrow\R^d$ is the proximal operator associated with the function $\g g$, and $\partial g(z)$ is the subifferential of $g$ at $z$~\cite{nesterov2003introductory}.
\end{prop}

\begin{remark}[{\sc Connection to Nesterov's smoothing technique}]
Alternatively, we can also view the MY envelope through the lens of Nesterov's smoothing technique~\cite{nesterov2005smooth}. More specifically, if Assumption~\ref{assumption:fg} is fulfilled, one can invoke a standard minimax theorem to verify that
\begin{equation*}
    g_\g(x) = \max_{z\in \R^d}\l\{ \langle x,z\rangle - 
    g^*(z) - \frac{\g}{2}\|z\|_2^2 \r\}, \qquad x\in \R^d,
\end{equation*}
where $g^*$ is the Fenchel conjugate of $g$. 
The right-hand side above plays a key role in Nesterov's technique for minimizing the non-smooth function~$F$ in \eqref{eq:pi}. 
\end{remark}

In view of the admissibility of~$F^{\MY}_{\gamma}$ by Proposition~\ref{prop:my}, applying ULA to the new potential~$F^{\MY}_{\gamma}$ leads to a well-defined algorithm, see Remark~\ref{rem:imp}. In addition, if~$\gamma$ is sufficiently small,~$\P^\MY_\g\propto e^{-{F_{\gamma}^\MY}}$  would be  close to the target distribution~$\P$ by Remark~\ref{rem:imp}.  This technique is known as MYULA~\cite{brosse2017sampling}.  However, a limitation of the MY envelope is that the minimizers of $F^{MY}_{\gamma}$ are not necessarily the same as the minimizers of~$F$. In turn,  the MAP estimator of~$\P_\g^\MY$, denoted by $x^\MY_\g$,  might not coincide with the MAP estimator of~$\P$, except in the limit of $\gamma\rightarrow 0$. That is 
\begin{equation*}
     \lim_{\g \rightarrow 0} F_{\gamma}(x_\g^{\MY})= \min_x F(x).\label{eq:onlyLimit}
\end{equation*}
This observation is particularly problematic because, as we will see later, very small values of $\g$ are often avoided in practice due to numerical stability issues.
In view of this discussion, our objective  is to replace the MY envelope with a new envelope that has the same minimizers as $F$ for \emph{all} sufficiently small $\g$, and not just in the limit of $\g\rightarrow0$.



\section{Forward-Backward Envelope}\label{sec:FB}

In this section, we will study an  envelope that addresses the limitations of the  MY envelope. More specifically, for $\g>0$, let us recall from \cite{stella2017forward} that the forward-backward (FB) envelope  of the function $F$ in~\eqref{eq:BigF} is defined as
\begin{equation}
    F_\g^{\FB}(x)  := f(x) -\frac{\g}{2}\|\nabla f(x)\|_2^2 + g_\g(x-\g\nabla f(x)), \qquad x\in \R^d,
    \label{eq:fbEnv}
    \tag{FB}
\end{equation}
where $g_\g$ was defined in \eqref{eq:MY}. A number of useful properties of $F^\FB_\g$ are collected below for the convenience of the reader~\cite{stella2017forward}. Recall that $P_{\g g}$ denotes the proximal operator associated with the function $\g g$ in~\eqref{eq:gradMY+proximal}. 
\begin{prop}[{\sc Properties of the FB envelope}]\label{prop:propsFB} Suppose that Assumption \ref{assumption:fg} is fulfilled. For $\g\in (0,1/\lambda_2)$ and every $x\in \R^d$, it holds that 
\begin{enumerate}[leftmargin=*,label={(\roman*)}]
    \item\label{item:sandwitch1} $F( P_{\g g}(x-\g \nabla f(x) ) \le F_\g^{\FB}(x) \le F(x) $, which relates the function $F$ to its FB envelope. 
    \item\label{item:sandwitch2} $F_{\frac{\g}{1-\g \lambda_2}}^{\MY}(x) \le  F_\g^{\FB}(x) \le F_\g^{\MY}(x)$, which relates the MY and FB  envelopes of the function $F$.
    \item\label{item:derFB} $F_\g^{\FB}$ is continuously differentiable and its gradient is given by $$\nabla F_\g^{\FB}(x)= \g^{-1}(I-\g \nabla^2 f(x)) (x - P_{\g g}(x-\g \nabla f(x))).$$ 
    \item\label{item:minsCoincide} $\argmin F^\FB_\g = \argmin F$, i.e., the function $F$ and its FB envelope have the same minimizers.  
\end{enumerate}

\end{prop}

In view of Proposition \ref{prop:propsFB}\ref{item:minsCoincide}, a remarkable property of the FB envelope is that the modes of $\P_\g^\FB \propto e^{-F^\FB_\g}$ coincide with the modes of the target measure $\P\propto e^{-F}$, for \emph{all} sufficiently small $\g$, rather than only in the limit of $\g\rightarrow 0$.  Indeed, very small values of $\g$ are often avoided in practice due to numerical stability issues. This observation signifies the advantage of the FB envelope over the MY envelope. Recall that the modes of the MY envelope coincide with those of $\P$ only in the limit of $\g\rightarrow 0$, see Section~\ref{sec:my-Env}. 

As a side note, let us remark that the proximal gradient descent algorithm for minimizing the (non-smooth) function $F$ coincides with the gradient descent (with variable metric) for minimizing the (smooth) function~$F^\FB_\g$, whenever $\g$ is sufficiently small~\cite{stella2017forward}.
It is also easy to use Proposition~\ref{prop:propsFB} to check the admissibility of the FB envelopes, as summarized below.  

\begin{prop}[{\sc Admissibility of FB envelopes}]\label{prop:fwbwAdmissible}  
Suppose that Assumption \ref{assumption:fg} is fulfilled. Then $\{F_\g^{\FB}:\g \in (0,\g^{\FB})\}$ are admissible envelopes of $F$ in~\eqref{eq:BigF}, where 
\begin{equation}
\g^{\FB}:=\frac{1}{2 \lambda_2+ 2\lambda_3(\lambda_0+R)}.
\label{eq:g0FB}
\end{equation}
Moreover, it holds that 
\begin{subequations}
\begin{align}
    & \l\| \nabla F^\FB_\g(x) - \nabla F^\FB_\g(y) \r\|_2
    \le \lambda_\g^\FB \|x-y\|_2,
\qquad x,y\in \R^d,  \label{eq:smoothLongStrCvx1} \\
    & \langle x-y,\nabla F^\FB_\g(x)-\nabla F^\FB_\g(y) \rangle \ge  \mu^\FB_\g \|x-y\|_2^2, \qquad \|x-y\|_2 \ge \rho_\g^\FB,
    \label{eq:smoothLongStrCvx2}
\end{align}\label{eq:smoothLongStrCvx}
\end{subequations}where 

\begin{align*}
    \lambda^{\FB}_\g &:= \g^{-1}+2\lambda_2+\lambda_3(\lambda_0+R),
    \quad \mu^\FB_\g:=  \lambda_2+ \lambda_3(\lambda_0+R), 
\\
\rho_\g^\FB &:= \frac{2R}{1-2\g( \lambda_2+ \lambda_3(\lambda_0+R)))}.
\end{align*}

\end{prop}

The equation \eqref{eq:smoothLongStrCvx} provides valuable information about the landscape of the FB envelope of $F$, which we now summarize: \eqref{eq:smoothLongStrCvx1} means that $F^\FB_\g$ is a $\lambda_\g^\FB$-smooth function. 
The smoothness of $F^\FB_\g$ in~\eqref{eq:fbEnv} is not surprising since both~$f$ and $g_\g$ are smooth functions. (Recall that  $g_\g$ is the MY envelope of $g$, which is known to be $\g^{-1}$-smooth.) 

Moreover, even though $F^\FB_\g$ is not necessarily a strongly convex function, \eqref{eq:smoothLongStrCvx2} implies that $F^\FB_\g$ behaves like a strongly convex function over long distances. As detailed in the proof,~\eqref{eq:smoothLongStrCvx2} holds essentially because the MY envelope of the indicator function $1_\CC$ is the function $\frac{1}{2\g}\dist(\cdot,\CC)^2$. The latter function grows quadratically faraway from the origin. Here, $\dist(\cdot,\CC)$ is the distance to the set~$\CC$. 

It is worth noting that a similar result to Proposition~\ref{prop:fwbwAdmissible} is implicit in~\cite{brosse2017sampling}. That is, the MY envelope~$F^\MY_\g$ also satisfies~\eqref{eq:smoothLongStrCvx}, albeit with different constants.  

\begin{remark}[{\sc Convergence in the Wasserstein metric}]\label{rem:cvgWass}
Recall from Remark~\ref{rem:imp} that ~$\P_\g^\FB$  converges weakly to $\P$ in the limit of $\g\rightarrow 0$. This weak convergence implies convergence in the Wasserstein metric by~\cite[Lemma~2.6]{eberle}: 
\begin{equation}
\lim_{\g\rightarrow 0} \W_1(\P_\g^\FB,\P)=0.
\label{eq:wassCvgWeak}
\end{equation}
We recall that, for two probability measures $\Q_1$ and $\Q_2$  that satisfy $\E_{x\sim \Q_1}\|x\|_2<\infty$ and $\E_{y\sim \Q_2}\|y\|_2 <\infty$, their $1$-Wasserstein or Kantorovich distance~\cite{villani2009optimal} is defined as 
\begin{align}
    \W_1(\Q_1,\Q_2):= \inf_{\substack{x\sim \Q_1\\y\sim \Q_2}} \E\|x-y\|_2.
    \label{eq:wassDist}
\end{align}
With some abuse of notation, throughout this work, we will occasionally replace the probability measures with the corresponding probability distributions or random variables. 
\end{remark}

A non-asymptotic version of Remark~\ref{rem:cvgWass} is presented below,  which bounds the Wasserstein distance between the two probability measures $\P^\FB_\g$ and $\P$. In effect, the result below is  an analogue of \cite[Proposition 5]{brosse2017sampling} for the MY envelope. The key ingredient of their result is the Steiner's formula for the volume of the set $\CC+\B(0,t)=\{x: \dist(x,\CC)\le t\}$ for every $t\ge 0$. The previous sum is in the Minkowski sense.  
Essentially, our proof strategy is to use Proposition~\ref{prop:propsFB}\ref{item:sandwitch2} to relate the FB and MY envelopes and then invoke \cite[Proposition 5]{brosse2017sampling}. 

\begin{thm}[{\sc Wassenstein distance between $\P^\FB_\g$ and $\P$}]\label{prop:distMeasures} Suppose that Assumption \ref{assumption:fg} is fulfilled.
For $\g\in (0,\g^{\FB})$,  it holds that 
\begin{align}
    \W_1(\P_\g^\FB,\P) & \le 
  \Cl[fcn]{dist1}+\frac{ \Cl[fcn]{dist3} R I_1(\g)+
  \Cr{dist3} I_2(\g)+
  \Cl[fcn]{dist2} I_2(\g/(1-\g\lambda_2))
  }{\vol(\CC)+I_1(\g)}+   \Cl[fcn]{dist4} R
  \nonumber\\
  & =: \Cl[fcn]{env}^\FB(\g)
  ,
  \label{eq:distProp}
\end{align}
where 
\begin{align}
    & \Cr{dist1}:=  \frac{e^{2\max_{x\in \CC} g(x)-\g\lambda_2 \min_z f(z)} \int_{\CC} \|x\|_2 e^{-(1-\g\lambda_2)f(x)} \, \der x }{\int e^{-f(x)-\frac{\dist(x,\CC)^2}{2\g}}\,\der x}\nonumber\\
    & \quad - \frac{e^{\g\lambda_2 \min_z f(z)-2\max_{x\in \CC} g(x)} \int_\CC \|x\|_2 e^{-f(x)}\, \der x}{\int e^{-(1-\g\lambda_2)f(x) - \frac{\dist(x,\CC)^2}{2\g/(1-\g\lambda_2)}}\, \der x},\nonumber\\
    & \Cr{dist3} := e^{\max_x f(x)-\min_x f(x)}, \nonumber\\
        & \Cr{dist2}:= e^{\max_x f(x) -  \min_x f(x) +2\max_{x\in \CC}g(x)}, \nonumber\\
    & \Cr{dist4}:= e^{\max_{x\in \CC}g(x)-\min_{x\in \CC}g(x) }-e^{\min_{x\in \CC}g(x)-\max_{x\in \CC} g(x)}, \nonumber\\
    & I_1(\g) :=  \sum_{i=0}^{d-1} \vol_{i}(\CC) \cdot (2\pi \g)^{\frac{d-i}{2}}, \nonumber\\
    & I_2(\g):= \sum_{i=0}^{d-1} \vol_i(\CC)\cdot  (2\pi\g)^{\frac{d-i}{2}} \l( \sqrt{\g(d-i+3)}+R \r).
    \label{eq:allConstantsDist}
\end{align}
Above, $\vol_i(\CC)$ is the $i$-th intrinsic volume of $\CC$, see~\cite{klain1997introduction}. In particular, the $d$-th volume of $\CC$ coincides with the standard volume of $\CC$, i.e.,  $\vol_d(\CC)=\vol(\CC)$. Moreover, to keep the notation light, above we have suppressed the dependence of $\Cr{dist1}$ to $\Cr{env}^{\FB}$ on $\CC,f,g,\g$. 
\end{thm}
As a sanity check, consider the special case of $g=1_\CC$, where $1_\CC$ is the indicator function for the set~$\CC$.   
Then we can use \eqref{eq:allConstantsDist} to verify that $\Cr{dist1}$  and $I_1(\g)$ and $I_2(\g)$ and $I_2(\g/(1-\g\lambda_2))$ all vanish when we send $\gamma\rightarrow 0$. Consequently, 
both the left- and right-hand sides of \eqref{eq:distProp} vanish if we send $\g\rightarrow0$. When $g=1_\CC$, then Theorem~\ref{prop:distMeasures} is precisely the analogue of  \cite[Proposition 5]{brosse2017sampling}. Their work, however, does not cover the case of~$g\ne 1_\CC$. 

In our result, when $g\ne 1_\CC$ and $\g\rightarrow0$, the right-hand side of \eqref{eq:distProp} converges to the nonzero value
\begin{align*}
    & \l(e^{2\max_{x\in \CC} g(x)} - e^{-2\max_{x\in \CC} g(x)}\r) \frac{\int_\CC \|x\|_2e^{-f(x)} \, \der x }{\int e^{-f(x)}\, \der x} \nonumber\\
    & \quad + \l( \frac{e^{\max_{x\in \CC} g(x)}}{e^{\min_{x\in \CC} g(x)}}
    -\frac{e^{\min_{x\in \CC} g(x)}}{ e^{\max_{x\in \CC} g(x)}}\r),
\end{align*}
unlike the left-hand side of~\eqref{eq:distProp}, which converges to zero by~\eqref{eq:wassCvgWeak}. Improving~\eqref{eq:distProp} in the case $g\ne 1_\CC$ appears to require highly restrictive assumptions on $g$ which we wish to avoid here. Moreover, in practice, {very} small values of~$\g$ are often avoided due to numerical stability issues. In this sense, improving \eqref{eq:distProp} for very small values of $\g$  might have limited practical value.

To summarize this section,  the FB envelopes $\{F^\FB_\g: \g\in (0,\g^{0,\FB})\}$ are admissible and we can use them as a differentiable proxy for the non-smooth function $F$ in \eqref{eq:BigF}. Crucially, the FB envelope addresses the key limitation of the MY envelope, i.e., the modes of  $\P^\FB_\g\propto e^{-F^\FB_\g}$ coincide with the modes $\P\propto e^{-F}$, for all sufficiently small $\g$, rather than only in the limit of $\g\rightarrow 0$. 


\section{EULA:\\Envelope Unadjusted Langevin Algorithm} \label{sec:newALg}

We have so far introduced two smooth envelopes for the non-smooth function~$F$ in \eqref{eq:BigF}, namely, the MY envelope $F^\MY_\g$ in \eqref{eq:MYEnv} and the FB envelope $F^\FB_\g$ in~\eqref{eq:fbEnv}. We also described in Section \ref{sec:FB} the advantage of the FB envelope over the MY envelope. To keep our discussion general, below we consider  admissible envelopes $\{F_\g:\g\in (0,\g^0)\}$ for the target function $F$ in~\eqref{eq:BigF}, see Definition~\ref{assump:envAssumption}. Our discussion below can be specialized to either of the envelopes by setting $F_\g=F^\MY_\g$ or $F_\g=F^\FB_\g$.


For the time being, let us fix $\g\in (0,\g^0)$. Unlike $F$, note that $\nabla F_\g$ exists and is Lipschitz continuous by Definition~\ref{assump:envAssumption}\ref{item:gradLips}. We can now use the ULA~\cite{DK19} to sample from $\P_\g\propto e^{-F_\g}$, as a proxy for the target measure $\P\propto e^{-F}$. 
The $k$-th iteration of the resulting algorithm is
\begin{equation}
    x_{k+1} = x_k - h \nabla F_{\g}(x_k)+ \sqrt{2h} \zeta_{k+1},
    \label{eq:vanillaULA}
\end{equation}
where $h$ is the step size and $\zeta_{k+1} \in \R^d$ is a standard  Gaussian random vector, independent of $\{\zeta_i\}_{i\le k}$. 
In particular, if we choose $F_\g=F^\MY_\g$, then \eqref{eq:vanillaULA} coincides with the MYULA from~\cite{brosse2017sampling}.

Under standard assumptions, to be reviewed later, the Markov chain $\{x_k\}_{k\ge 0}$ in \eqref{eq:vanillaULA} has a unique invariant probability measure, which we denote by $\wh{\P}_{\g,h}$. There are two sources of error that contribute to the difference between $\wh{\P}_{\g,h}$ and  the target measure $\P$ in~\eqref{eq:pi}, which we list below:
\begin{enumerate}[leftmargin=*]
    \item First, note that \eqref{eq:vanillaULA} is only intended to sample from $\P_\g$, as a proxy for the target distribution~$\P$. That is, the first source of error is the difference between the two probability measures $\P_\g$ and $\P$. 
    \item Second, the step size $h$ is known to contribute to the difference between the two probability measures $\wh{\P}_{\g,h}$ and~$\P_\g$, see~\cite{DK19}. This bias  vanishes only in the limit of $h\rightarrow 0$. 
\end{enumerate}
In fact, instead of \eqref{eq:vanillaULA}, we study here a slightly more general algorithm that allows $\g$ and $h$ to vary.
More specifically, for a nonincreasing sequence $\{\g_k\}_{k\ge 0}$ and step sizes $\{h_k\}_{k\ge 0}$, the $k$-th iteration of this more general algorithm is 
\begin{align}
    x_{k+1} = x_k - h_k \nabla F_{\g_k}(x_k) + \sqrt{2h_k} \zeta_{k+1},
    \label{eq:eula}
    \tag{EULA}
\end{align}
where $\zeta_{k+1}in \R^d$ is a standard Gaussian random vector independent of $\{\zeta_i\}_{i\le k}$. 
\eqref{eq:eula} stands for Envelope Unadjusted Langevin Algorithm. 

In particular, if we set $\g_k=\g$ in~\eqref{eq:eula} for every $k\ge0$, then we retrieve~\eqref{eq:vanillaULA}.
Alternatively, if $\{\g_k\}_{k\ge 0}$ is a decreasing sequence, then $F_{\g_k}$ becomes an increasingly better approximation of the target potential function~$F$ as $k$ increases, see Definition~\ref{assump:envAssumption}\ref{item:pntCvg}. That is, \eqref{eq:eula} uses increasingly better approximations of the potential function $F$ as $k$ increases. 

We next present the iteration complexity of the \eqref{eq:eula} for admissible envelopes $\{F_\g:\g\in (0,\g^0)\}$, where admissibility was defined in Definition~\ref{assump:envAssumption}.  
The result below can be specialized to both MY and FB envelopes by setting $F_\g=F^\MY_\g$ or $F_\g=F^\FB_\g$, respectively.  
\begin{thm}[{\sc Iteration complexity of \eqref{eq:eula}}]\label{thm:eulaComp}
For $\g^0>0$, consider admissible envelopes $\{F_\g:\g\in (0,\g^0)\}$ of $F$ in~\eqref{eq:BigF}, see Definition~\ref{assump:envAssumption}.  
For $\mu_\g>0$ and $\rho_\g\ge 0$, we additionally assume that $F_\g$ satisfies the inequality 
\begin{equation}
    \langle x-y,\nabla F_\g(x)-\nabla F_\g(y) \rangle \ge \mu_\g \|x-y\|_2^2, \qquad \|x-y\|_2\ge \rho_\g, \qquad \g\in (0,\g^0).
    \label{eq:longDistanceStrCvx}
\end{equation}
Consider two sequences $\{\g_k\}_{k\ge 0}\subset (0,\g^0)$ and $\{h_k\}_{k\ge 0} \subset \R_+$. For the algorithm \eqref{eq:eula}, let $\Q_k$ denote the law of $x_k$ for every integer $k\ge 0$. That is, $x_k \sim \Q_k$ for every $k\ge 0$. Then the $\W_1$ distance between $\Q_k$ and the target measure $\P\propto e^{-F}$ in~\eqref{eq:pi} is bounded by 
\begin{align}
    \W_1(\Q_{k},\P) 
 & \le  e^{\Cl[fcn]{radius}\Cl[fcn]{q}} \prod_{i=0}^{k-1} (1-\Cl[fcn]{contraction} h_i) \cdot  \W_1(\Q_0,\P_{\g_{0}}) + 
    e^{\Cr{radius}\Cr{q}}\sum_{i=0}^{k-1} \a_i \prod_{j=i+1}^{k-1}(1-\Cr{contraction}h_j) \nonumber \\
    &+ \W_1(\P_{\g_k},\P),
\end{align}
for every $k\ge 0$, provided that 
\begin{align*}
\g_k \in (0,\g^0), \quad  
     h_k \le \frac{1}{\lambda_{\g_k}} \min\l( \frac{1}{6}, \frac{\mu_{\g_k}}{\lambda_{\g_k}} , \frac{\lambda_{\g_k} \rho_{\g_k}^2}{3}, \frac{c_0^2}{970 \lambda_{\g_k} \rho_{\g_k}^2}\r), \quad  k\ge 0.
\end{align*}
Above, $c_0\ge 0.007$ is a universal constant specified in \cite[Equation (6.6)]{eberle2019quantitative}. Moreover, 
\begin{align*}
 & \Cr{radius} := (1+h_k\lambda_{\g_{k}})\rho_{\g_k} \le 7 \rho_{\g_k}/6, 
    \quad \Cr{q}:= 7 \lambda_{\g_k}\rho_{\g_k}/c_0, \nonumber\\
    & \Cr{contraction} := \min\l( \frac{\mu_{\g_k}}{2},\frac{245}{24c_0} (\lambda_{\g_k} \rho_{\g_k})^2 \r) e^{-\frac{49}{6c_0} \lambda_{\g_k} \rho_{\g_k}^2}, \nonumber\\
 & \a_k := \lambda_{\g_{k}} \sqrt{h_k^3 d} \cdot \l( \sqrt{h_k \lambda_{\g_{k}} } +\sqrt{2}\r)
    +\W_1(\P_{\g_k},\P)+\W_1(\P_{\g_{k+1}},\P), \,\,\, k\ge 0.
\end{align*}

\end{thm}
Note that, when $\rho_\g=0$, then \eqref{eq:longDistanceStrCvx} requires $F_\g$ to be $\mu_\g$-strongly convex for every $\g\in (0,\g^0)$. This can happen, for example, when $f$ itself is a strongly convex function. A particularly important special case of Theorem~\ref{thm:eulaComp} is when we choose the FB envelope, and use a fixed $\g$ and step size $h$.

\begin{corr}[{\sc Iteration complexity of \eqref{eq:eula} for FB envelope}]\label{cor:eula}Suppose that Assumption~\ref{assumption:fg} is fulfilled. For the algorithm \eqref{eq:eula}, suppose that $\g_k=\g\in (0,\g^\FB)$ and  $F_{\g_k}=F^\FB_\g$ and $h_k=h>0$ for every integer $k\ge 0$, see~\eqref{eq:fbEnv} and~\eqref{eq:g0FB}. In~\eqref{eq:eula}, also let $\Q_k$ denote the law of $x_k$ for every $k\ge 0$. That is, $x_k \sim \Q_k$ for every integer $k\ge 0$. Then the $\W_1$ distance between $\Q_k$ and the target measure $\P\propto e^{-F}$ in~\eqref{eq:pi} is bounded by 
\begin{align}
    \W_1(\Q_{k},\P) 
 & \le  e^{\Cr{radius}^\FB\Cr{q}^\FB}  (1-\Cr{contraction}^\FB h)^k \cdot  \W_1(\Q_0,\P_{\g}) +
    \frac{\a^\FB e^{\Cr{radius}^\FB\Cr{q}^\FB}}{\Cr{contraction}^\FB h}+ \Cr{env}^\FB,
\end{align}
for every $k\ge 0$, provided that 
\begin{align}
\g\in (0,\g^0), \qquad 
     h \le \frac{1}{\lambda_{\g}^\FB} \min\l( \frac{1}{6}, \frac{\mu_{\g}^\FB}{\lambda_{\g}^\FB} , \frac{\lambda_{\g}^\FB (\rho_{\g}^\FB)^2}{3}, \frac{c_0^2}{970 \lambda_{\g}^\FB (\rho_{\g}^
     \FB)^2}\r).
\end{align}
Above, $c_0\ge 0.007$ is a universal constant specified in \cite[Equation (6.6)]{eberle2019quantitative}. Moreover, 
\begin{align*}
 & \Cr{radius}^\FB := (1+h\lambda_{\g}^\FB)\rho_{\g}^\FB \le 7 \rho_{\g}^\FB/6, 
    \quad \Cr{q}^\FB:= 7 \lambda_{\g}^\FB\rho_{\g}^\FB/c_0, \nonumber\\
    & \Cr{contraction}^\FB:= \min\l( \frac{\mu_{\g}^\FB}{2},\frac{245}{24c_0} (\lambda_{\g}^\FB \rho_{\g}^\FB)^2 \r) e^{-\frac{49}{6c_0} \lambda_{\g}^\FB (\rho_{\g}^\FB)^2}, \nonumber\\
 & \a^\FB := \lambda_{\g}^\FB \sqrt{h^3 d} \cdot \l( \sqrt{h \lambda_{\g }^\FB} +\sqrt{2}\r)
    +2\Cr{env}^\FB.
\end{align*}
The remaining quantities were defined in Propositions \ref{prop:fwbwAdmissible} and \ref{prop:distMeasures}. 

\end{corr}

We remark that Corollary~\ref{cor:eula} for the FB envelope is the analogue of \cite[Proposition~7]{brosse2017sampling} for the MY  envelope. However, note that \cite[Proposition~7]{brosse2017sampling} requires $f$ to be strongly convex whereas we merely assume $f$ to be convex, see Assumption~\ref{assumption:fg}. 

\section{Proof of Theorem \ref{thm:eulaComp}\\ (Iteration Complexity of EULA)}\label{sec:proofMain}

To begin, we let $Q_k$ denote the Markov transition kernel associated with the Markov chain $\{x_k\}_{k\ge 0}$. This transition kernel is specified as 
\begin{align}
    Q_k(x,\cdot) := \text{Normal}(x-h_k\nabla F_{\g_k}(x),2h_{k}I_d), \qquad x\in \R^d,
    \label{eq:QTransitionKernel}
\end{align}
where $\text{Normal}(a,B)$ is the Gaussian probability measure with mean $a\in \R^d$ and covariance matrix $B\in \R^{d\times d}$. {Above, note that $Q_k$ depends on both $h_k$ and $\gamma_k$.}
We also let $\Q_k$ denote the law of $x_k$,  i.e., $x_k \sim \Q_k$. Using the standard notation, we can now write that 
\begin{equation}
\Q_{k+1}=\Q_k Q_{{k}}, \qquad  k\ge 0.
\label{eq:Qq}
\end{equation}
To be precise, \eqref{eq:Qq} is equivalent to 
\begin{align}
    \Q_{k+1}(\der y) = \int \Q_{k}(\der x) Q_{k}( x,\der y) , \qquad y\in \R^d,\, k\ge 0.
\end{align}
Recall that $\P_{\g_{k+1}}$ serves as a proxy for the target probability measure $\P$.
The~$\W_1$ distance between $\Q_{k+1}$ and $\P_{\g_{k}}$ can be bounded as 
\begin{align}
    \W_1(\Q_{k+1},\P_{\g_{k}}) & = \W_1( \Q_{k}Q_{k},\P_{\g_{k}})  
    \qquad \text{(see \eqref{eq:Qq})}
    \nonumber\\
    & \le 
    \W_1( \Q_{\g_k}Q_{k},\P_{\g_{k}}Q_{k}) + \W_1( \P_{\g_{k}}Q_{k}, \P_{\g_{k}}),
    \label{eq:mainTri}
\end{align}
where the second line follows from the triangle inequality.
We separately control each $\W_1$ distance in the last line above. For the first distance, we plan to invoke Theorem 2.12 from \cite{eberle2019quantitative}, reviewed below for the convenience of the reader. It is worth noting that a similar result to the one below appears in \cite[Corollary 2.4]{majka2020nonasymptotic}.
\begin{prop}[{\sc \cite{eberle2019quantitative}, Theorem 2.12}]\label{prop:eberle} Let 
\begin{align*}
    & \Cr{radius} := (1+h_k\lambda_{\g_{k}})\rho_{\g_k} \le 7 \rho_{\g_k}/6, 
    \quad \Cr{q}:= 7 \lambda_{\g_k}\rho_{\g_k}/c_0, \nonumber\\
    & \theta(r) :=\int_0^r e^{-\Cr{q}\min(s,r_1)} \der s, \qquad \Theta(x,y):= \theta(\|x-y\|_2),\nonumber\\
    & \Cr{contraction} := \min\l( \frac{\mu_{\g_k}}{2},\frac{245}{24c_0} (\lambda_{\g_k} \rho_{\g_k})^2 \r) e^{-\frac{49}{6c_0} \lambda_{\g_k} \rho_{\g_k}^2}, \nonumber\\
    & h_k \le \frac{1}{\lambda_{\g_k}} \min\l( \frac{1}{6}, \frac{\mu_{\g_k}}{\lambda_{\g_k}} , \frac{\lambda_{\g_k} \rho_{\g_k}^2}{3}, \frac{c_0^2}{970 \lambda_{\g_k} \rho_{\g_k}^2}\r),
\end{align*}
where $c_0\ge 0.007$ is a universal constant specified in \cite[Equation (6.6)]{eberle2019quantitative}.
Then it holds that 
\begin{align}
    \W_\Theta(\Q_{k}Q_{k},\P_{\g_{k}}Q_{k}) \le (1-\Cr{contraction} h_k) \cdot  \W_\Theta(\Q_k,\P_{\g_{k}}), 
\end{align}
where $\W_\Theta$ is defined similar to \eqref{eq:wassDist} but the $\ell_2$-norm is replaced with $\Theta$. Above, to keep the notation light, we have suppressed the dependence of $\Cr{radius}$ and $\Cr{q}$ and $\Cr{contraction}$ on $\CC,f,g,h_k$.
Moreover, the two metrics $\W_\Theta$ and~$\W_1$ are related as 
\begin{equation}
    e^{-\Cr{radius}\Cr{q}} \W_1 \le \W_\Theta\le \W_1.
    \label{eq:W1WTheta}
\end{equation}

\end{prop}

For the second $\W_1$ distance in the last line of \eqref{eq:mainTri}, the following result is standard, {see appendix for the proof.} 
\begin{lem}[{\sc Discretization error}]\label{lem:termTwo}
It holds that 
\begin{equation}
    \W_1( \P_{\g_{k}}Q_{k}, \P_{\g_{k}}) \le 
    \Cl[fcn]{disc}
    :=\lambda_{\g_{k}} \sqrt{h_k^3 d} \cdot \l( \sqrt{h_k \lambda_{\g_{k}} } +\sqrt{2}\r). 
    \label{eq:discErr}
\end{equation}
\end{lem}
In fact, it is more common to write the left-hand side of \eqref{eq:discErr} in terms of the Markov transition kernel of the corresponding Langevin diffusion, as discussed in the proof of Lemma~\ref{lem:termTwo}. 
By combining Proposition~\ref{prop:eberle} and Lemma~\ref{lem:termTwo}, we can now revisit \eqref{eq:mainTri} and write that  
\begin{align}
     \W_\Theta(\Q_{k+1},\P_{\g_{k}}) 
    & \le 
    \W_\Theta( \Q_{k}Q_{k},\P_{\g_{k}}Q_{k}) + \W_\Theta( \P_{\g_{k}}Q_{k}, \P_{\g_{k}})
    \qquad \text{(see \eqref{eq:mainTri})} \nonumber\\
    & \le (1-\Cr{contraction} h_k)  \W_\Theta(\Q_k,\P_{\g_{k}}) + \W_1( \P_{\g_{k}}Q_{k}, \P_{\g_{k}}) 
    \qquad \text{(see \eqref{eq:W1WTheta})}
    \nonumber\\
    & \le (1-\Cr{contraction}h_k)  \W_\Theta(\Q_k,\P_{\g_{k}}) +\Cr{disc}.
    \qquad \text{(Lemma \ref{lem:termTwo})}
    \label{eq:OneIteration}
\end{align}
Using the triangle inequality, it immediately follows that 
\begin{align}
   & \W_\Theta(\Q_{k+1},\P_{\g_{k+1}}) \nonumber\\
   & \le \W_\Theta(\Q_{k+1},\P_{\g_k}) +\W_\Theta(\P_{\g_k},\P)+\W_\Theta(\P_{\g_{k+1}},\P) \qquad \text{(triangle inequality)} \nonumber\\
    & \le (1-\Cr{contraction} h_k)  \W_\Theta(\Q_k,\P_{\g_{k}}) + \Cr{disc}
    +\W_\Theta(\P_{\g_k},\P)+\W_\Theta(\P_{\g_{k+1}},\P) 
    \qquad \text{(see \eqref{eq:OneIteration})}
    \nonumber\\
    & \le (1-\Cr{contraction} h_k)  \W_\Theta(\Q_k,\P_{\g_{k}}) + \Cr{disc}
    +\W_1(\P_{\g_k},\P)+\W_1(\P_{\g_{k+1}},\P) \qquad \text{(see \eqref{eq:W1WTheta})} \nonumber\\
    & =:(1-\Cr{contraction} h_k)  \W_\Theta(\Q_k,\P_{\g_{k}}) +\a_k.  
    \label{eq:OneIteration2}
\end{align}
By unwrapping \eqref{eq:OneIteration2}, we find that 
\begin{align}
      \W_1 (\Q_{k},\P_{\g_{k}}) 
     & \le e^{\Cr{radius}\Cr{q}} \W_\Theta(\Q_{k},\P_{\g_{k}}) 
     \qquad \text{(see \eqref{eq:W1WTheta})}
     \nonumber\\
    & \le e^{\Cr{radius}\Cr{q}} \prod_{i=0}^{k-1} (1-\Cr{contraction}h_i) \cdot   \W_\Theta(\Q_0,\P_{\g_{0}}) \nonumber\\
& \qquad \quad     + 
    e^{\Cr{radius}\Cr{q}}\sum_{i=0}^{k-1} \a_i \prod_{j=i+1}^{k-1}(1-\Cr{contraction}h_j)
    \qquad \text{(see \eqref{eq:OneIteration2})}
    \nonumber\\
    & \le e^{\Cr{radius}\Cr{q}} \prod_{i=0}^{k-1} (1-\Cr{contraction}h_i) \cdot   \W_1(\Q_0,\P_{\g_{0}}) \nonumber\\
    & \qquad \quad + 
    e^{\Cr{radius}\Cr{q}}\sum_{i=0}^{k-1} \a_i \prod_{j=i+1}^{k-1}(1-\Cr{contraction}h_j).
    \qquad \text{(see \eqref{eq:W1WTheta})}
    \label{eq:penUlt}
\end{align}
Lastly, we can use \eqref{eq:penUlt} in order to bound the $\W_1$ distance at iteration $k$ to the target measure $\P$  as
\begin{align}
    \W_1(\Q_{k},\P) & \le \W_1(\Q_{k},\P_{\g_{k}}) + \W_1(\P_{\g_{k}}, \P) 
    \qquad \text{(triangle inequality)}
    \nonumber\\
 & \le  e^{\Cr{radius}\Cr{q}} \prod_{i=0}^{k-1} (1-\Cr{contraction}h_i) \cdot  \W_1(\Q_0,\P_{\g_{0}}) \nonumber\\ 
 & \qquad + 
    e^{\Cr{radius}\Cr{q}}\sum_{i=0}^{k-1} \a_i \prod_{j=i+1}^{k-1}(1-\Cr{contraction}h_j) + \W_1(\P_{\g_k},\P),
\end{align}
which completes the proof of Theorem \ref{thm:eulaComp}.

\section{Numerical experiments}\label{sec:numerics}

A number of numerical experiments are 
{presented below to support our theoretical contributions.}

\subsection{Truncated Gaussian}\label{sec:truncGauss}
Our first numerical experiment deals with sampling from a truncated Gaussian distribution, restricted to
{a box
$K_{d} \subset \R^d$.} For this problem the potential {$U:\R^d\rightarrow\R$} is {specified as}
\begin{equation}
    U(x):=\frac{1}{2} \left \langle  x,\Sigma^{-1}x \right \rangle+\iota_{K_{d}}(x).
\end{equation}
Here similarly to  \cite{brosse2017sampling}
 the {$(i,j)$th entry of the} covariance matrix is given by
\[
\Sigma_{i,j}:=\frac{1}{1+|i-j|}.
\]
We now
consider three scenarios, {namely,}  
\begin{itemize}
    \item   $d=2$  {with} $K_{2}=[0.5] \times [0,1]$,
    \item $d=10$ {with} $K_{10}=[0,5] \times [0,0.5]^{9}$
    \item  $d=100$ {with} $K_{100}=[0,5] \times [0,0.5]^{99}$.  
\end{itemize} 

\begin{figure}
    \centering
    \subfigure[$\mathbb{E}(x_{1})$]
    {
        \includegraphics[width=.45\textwidth]{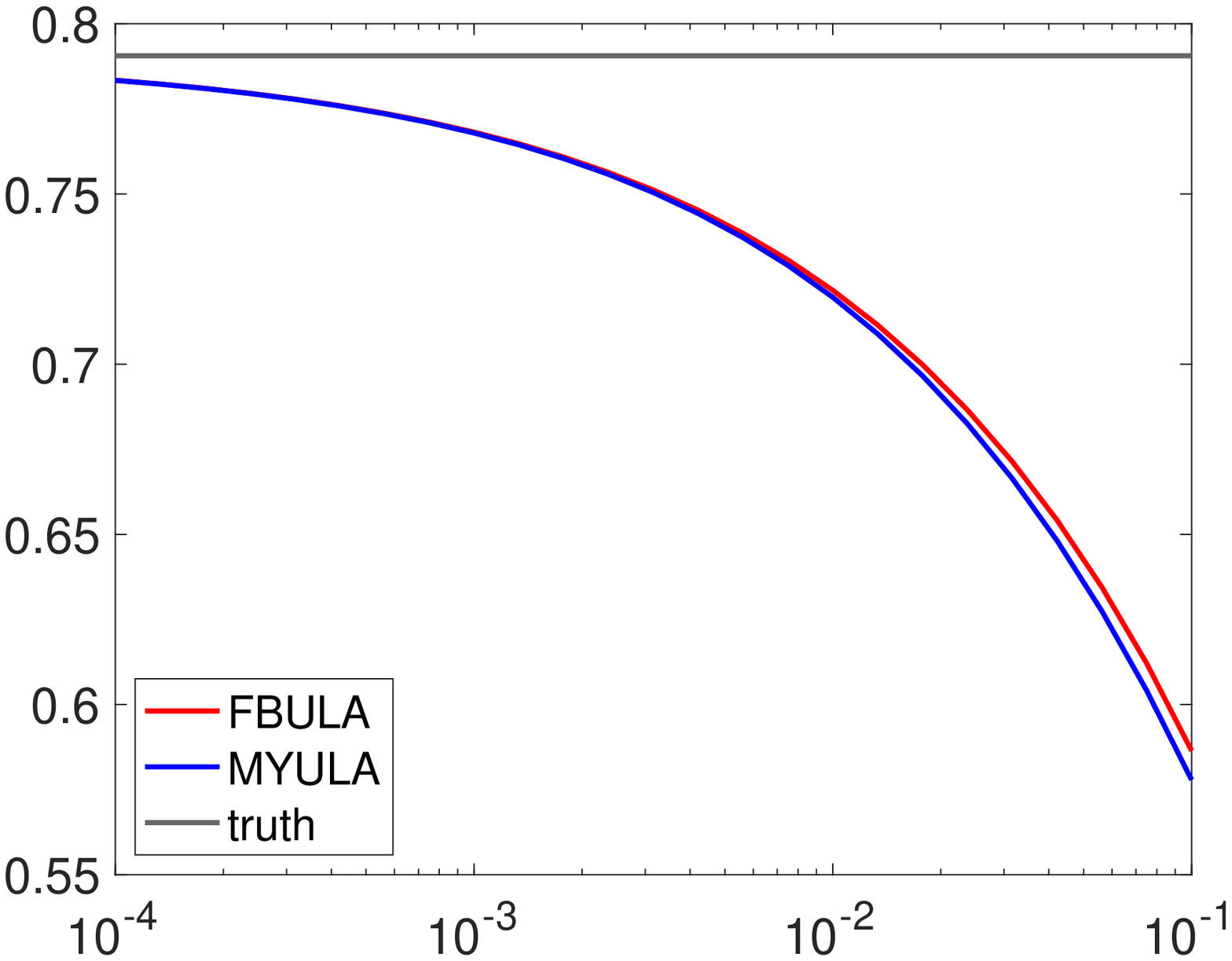}
        \label{subfig:gamma_mean}
    }\subfigure[$\text{Var}(x_{1})$]
    {
        \includegraphics[width=.45\textwidth]{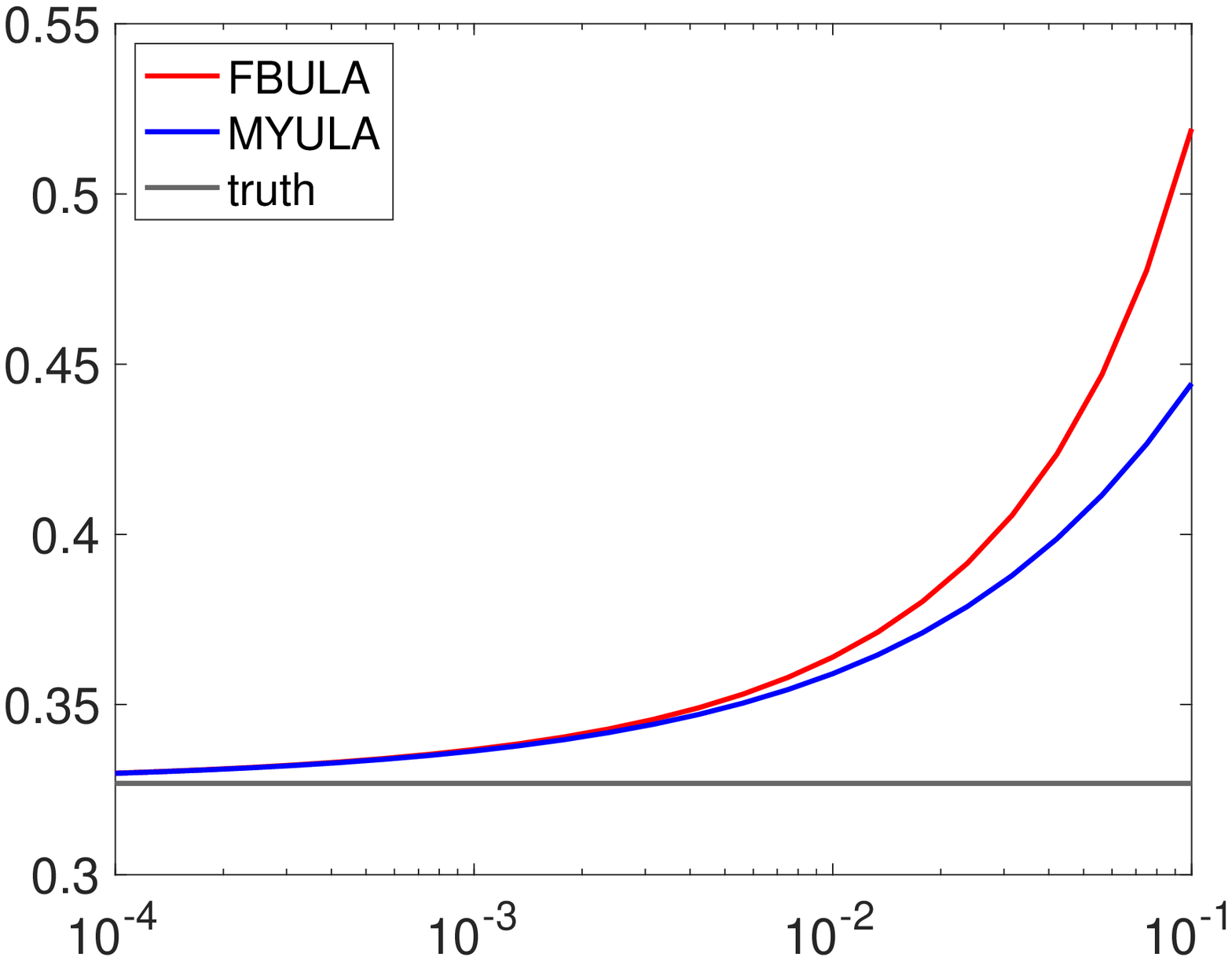}
        \label{subfig:gamma_var} 
    }\\
    \caption{
    {\small This figure compares the MY and FB envelopes for the two-dimensional truncated Gaussian distribution $\P\propto e^{-U-1_{K_2}}$ specified in Section~\ref{sec:truncGauss}. The horizontal lines in the left and right panels show, respectively, the expectation and variance of the first coordinate, namely, $\E_\P[x_1]$ and $\mathrm{var}_\P[x_1]=\E_\P[x_1^2]-(\E_\P[x_1])^2$. The blue and red graphs in both panels show the estimated values of $\E_\P[x_1]$ and $\mathrm{var}_\P[x_1]$, obtained via MY and FB envelopes. That is, the graphs on the left correspond to $\E_{\P_\gamma^\MY}[x_1]$ and $\E_{\P_\gamma^\FB}[x_1]$, for various values of $\gamma$. Similarly, the graphs on the right correspond to $\mathrm{var}_{\P_\gamma^\MY}[x_1]$ and $\mathrm{var}_{\P_\gamma^\FB}[x_1]$, for various values of $\gamma$.    } 
    }
    \label{fig:gamma_d2_comp}
\end{figure}

{Using quadrature techniques,} it is possible in the two-dimensional case {($d=2$)} to calculate exactly the mean and the covariance of the 
truncated {Gaussian} distribution, as well as the corresponding approximations {obtained via MY and FB envelopes.} {More specifically,} 
Figure~\ref{fig:gamma_d2_comp} {uses MATLAB's \texttt{integral2} command to plot the following quantities for various values of $\g$:
\begin{align*}
   \E_{\P^{\MY}_\gamma}[x_1] := \frac{ \int x_1 e^{-F^{\MY}_\g(x)} \,\der x}{\int e^{-F^{\MY}_\g(x)}\,\der x},
\quad
 \mathrm{var}_{\P^{\MY}_\gamma}[x_1]:= \frac{ \int x_1^2 e^{-F^{\MY}_\g(x)} \,\der x}{\int e^{-F^{\MY}_\g(x)}\,\der x}-  (\E_{\P^{\MY}_\gamma}[x_1] )^2,
 \\
   \E_{\P^{\FB}_\gamma}[x_1] := \frac{ \int x e^{-F^{\FB}_\g(x)} \,\der x}{\int e^{-F^{\FB}_\g(x)}\,\der x},
   \quad
   \mathrm{var}_{\P^{\FB}_\gamma}[x_1]:= \frac{ \int x_1^2 e^{-F^{\FB}_\g(x)} \,\der x}{\int e^{-F^{\FB}_\g(x)}\,\der x}- (\E_{\P^{\FB}_\gamma}[x_1] )^2.
\end{align*}
The horizontal lines in Figure~\ref{fig:gamma_d2_comp} show the ground truth values obtained by MATLAB's \texttt{integral2} command, i.e., 
\begin{align*}
    \E_{\P}[x_1] := \frac{ \int x_1 e^{-F(x)} \,\der x}{\int e^{-F(x)}\,\der x},
    \quad
    \mathrm{var}_{\P}[x_1]:= \frac{ \int x_1^2 e^{-F(x)} \,\der x}{\int e^{-F(x)}\,\der x}-  (\E_{\P}[x_1] )^2.
\end{align*}}For small values of the parameter 
$\gamma$, {we observe in Figure~\ref{fig:gamma_d2_comp} that}
{the FB envelope  better approximates} the  mean of the first component  than 
{the MY envelope.}
{However, the FB envelope tends} to overestimate the variance. This can be understood by comparing the two envelopes since in the case of the MY envelope the smoothing is more localized compared to the FB envelope.


\begin{figure}[h!]
    \centering
    \includegraphics[scale=.5]{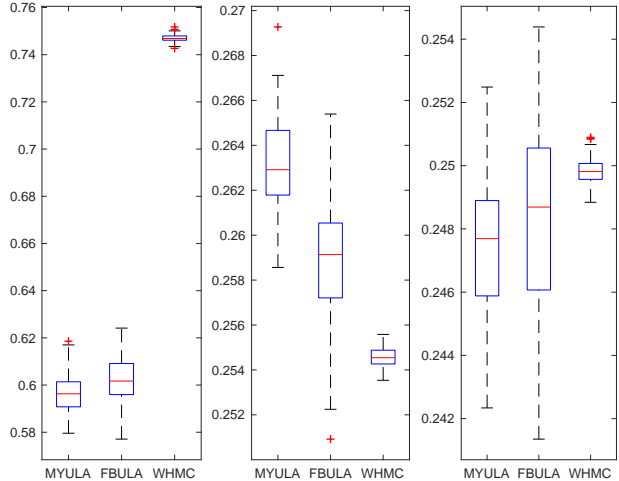}
    \caption{\small This figure shows the boxplots for {the expectations of} $x_{1},x_{2},x_{3}$ for the truncated Gaussian {distribution} in dimension~$10$ obtained by MYULA, FBULA, and wHMC. The last approach serves as the ground truth.}
    \label{fig:boxplot_d10}
\end{figure}

\begin{figure}
    \centering
    \includegraphics[scale=.5]{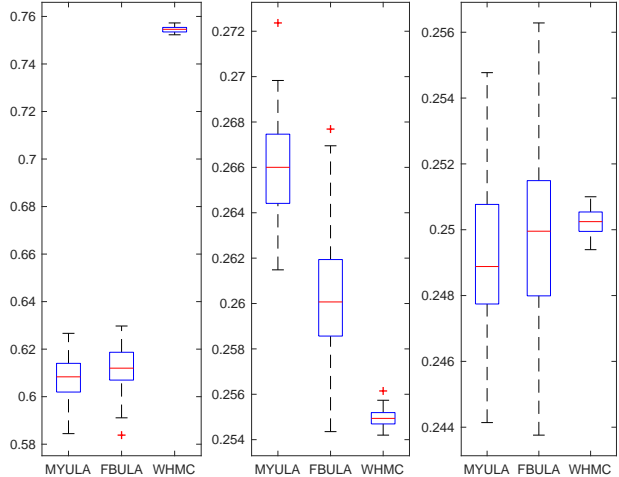}
    \caption{{\small
    Boxplots for {the expectations of} $x_{1},x_{2},x_{3}$ for the truncated Gaussian {distribution} in dimension~$100$, { obtained by MYULA, FBULA, and wHMC. The last approach serves as the ground truth.}} 
    }
    \label{fig:boxplot_d100}
\end{figure}

Such explicit calculations are not 
{tractable in higher dimensions, i.e., for $d\in \{10,100\}$.}
{Instead,} we now generate $10^6$ samples from the truncated Gaussian distribution $\P$ by applying MYULA and FBULA. 
{As our ground truth,} we also generate $10^{5}$ samples {from $\P$} with  the wall HMC (wHMC)~\cite{PP14}.  
{In all three approaches,} the initial $10\%$ {of the obtained samples} 
{are}
discarded as {the} burn-in period. In  terms of the parameters, we set $\gamma=0.05 $ and fix $h=0.005$
for all {of} our experiments.

{The results are visualized  in Figures~\ref{fig:boxplot_d10} and \ref{fig:boxplot_d100}. More specifically, Figure~\ref{fig:boxplot_d10} corresponds to $d=10$ and shows the estimates for $\E_\P[x_i]$~for $i\in \{1,2,3\}$, obtained by MYULA, FBULA, and wHMC. Similarly, Figure~\ref{fig:boxplot_d100} corresponds to $d=100$.}
{These figures indicate that,}
{in all of these} cases, FBULA is providing a more accurate approximation of the expectation {compared to} MYULA.

\subsection{Tomographic image reconstruction}

\begin{figure}
    \centering
    \subfigure[True image $x$]
    {
        \includegraphics[scale=.25]{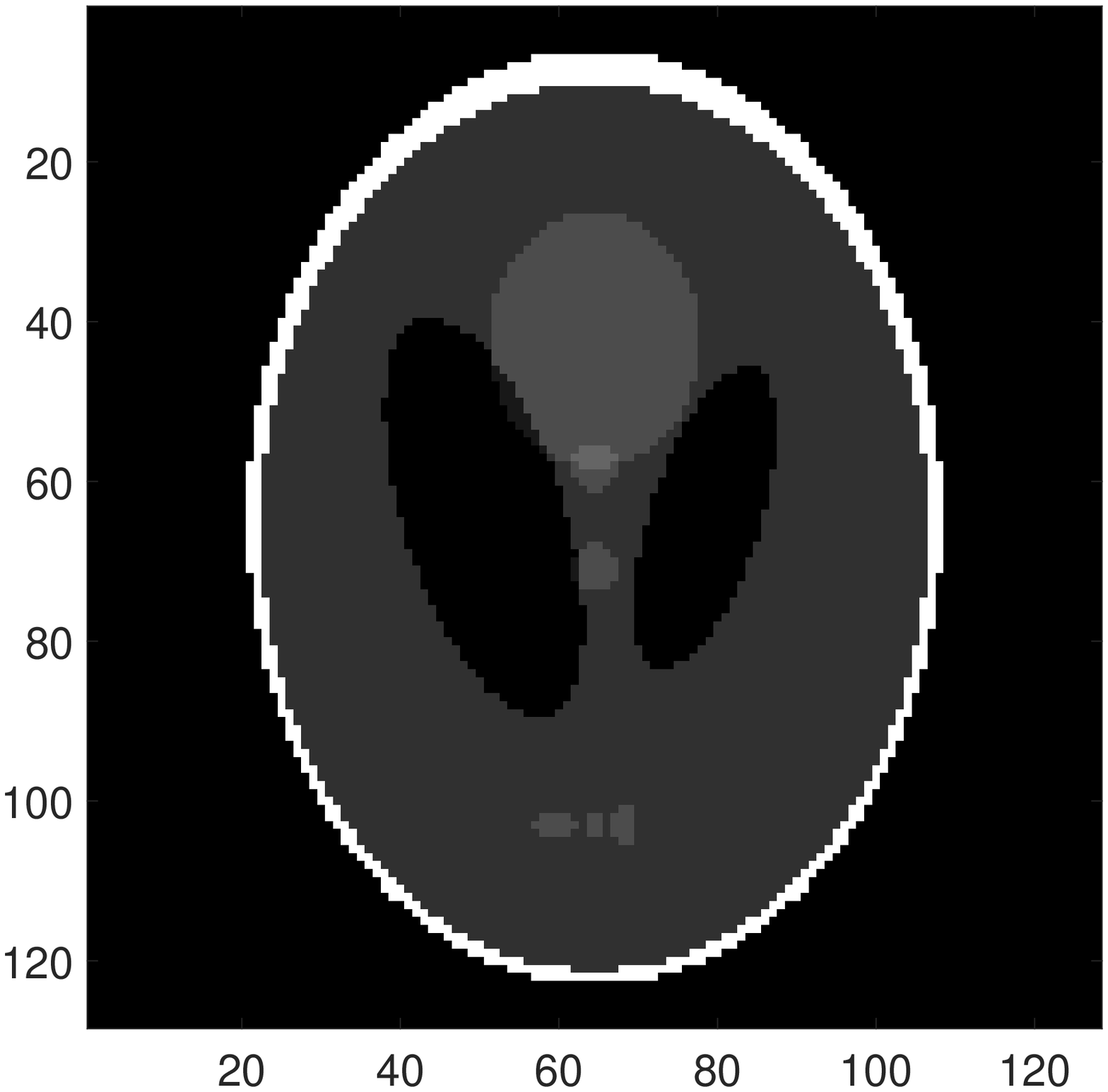}
        \label{subfig:tomography_x}
    }
    \subfigure[Observation $y$]
    {
        \includegraphics[scale=.25]{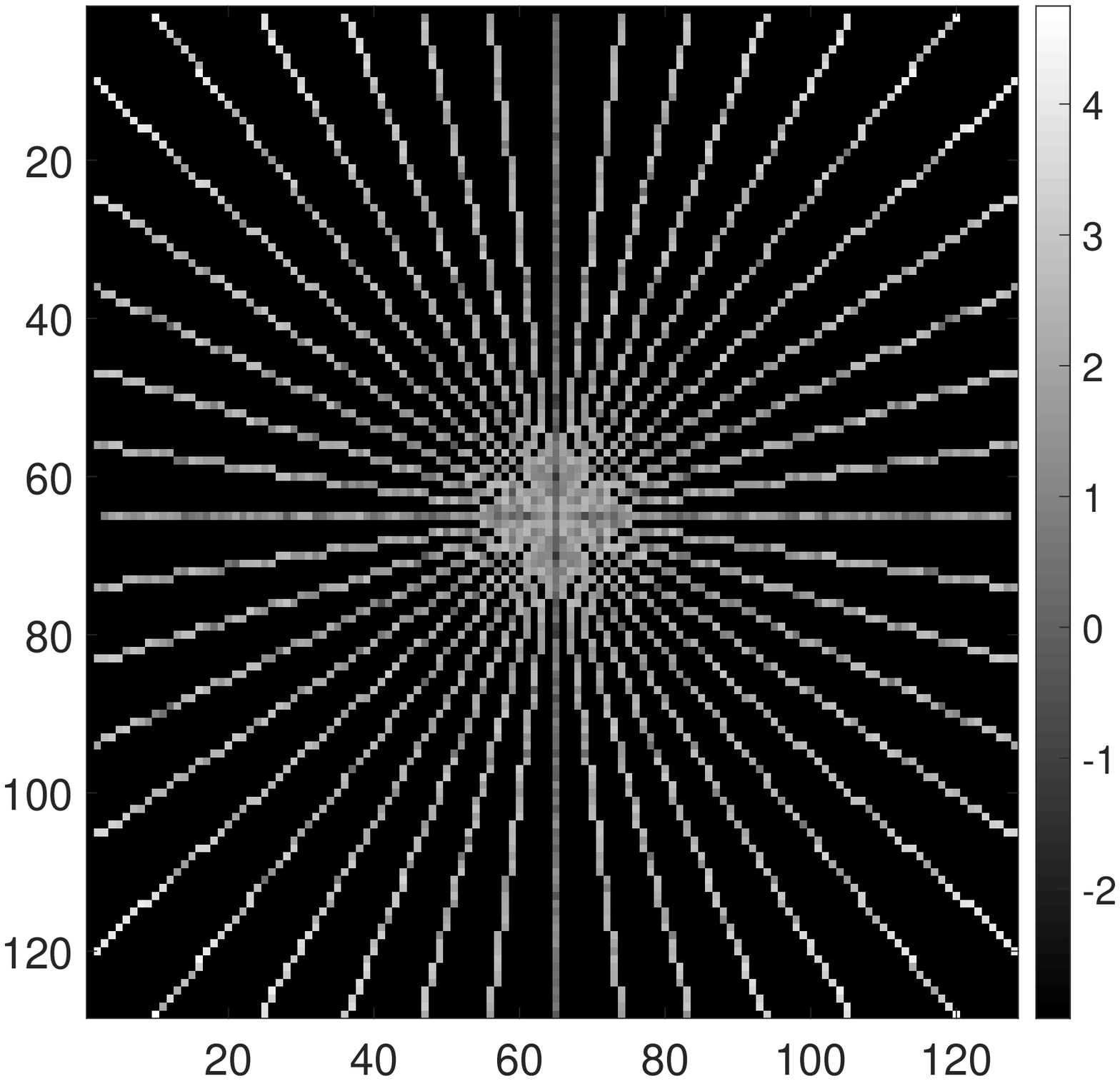}
        \label{subfig:tomography_y_amp_four_coeff} 
    }
    \caption{\textit{Tomography} experiment: (a) True image $x$ of dimension $d=128\times 128$. (b) Incomplete and noisy observation $y$, amplitude of Fourier coefficients in logarithmic scale.}
    \label{fig:tomography_x_y}
\end{figure}

We now study a tomographic image reconstruction problem. In this case the true image is taken to  the \textit{Shepp-Logan phantom} test image of dimension $d=128 \times 128$, in which we applied a Fourier operator~$F$ followed by a subsampling operator $A$, reducing the observed pixels 
{by}
approximately $85\%$. Finally, zero-mean additive Gaussian noise $\xi$ is added with standard deviation $\sigma = 10^{-2}$ to produce an incomplete observation $y=AFx + \xi$ where $y \in \mathbb{C}^p$. Note that $p < d$. With regards to the prior, we use the total-variation norm with an additional constraint for the size of the pixels. This leads to  the following posterior distribution:
\begin{equation} \label{eq:tomography_posterior_dist}
    \pi (x) \propto \exp \left[ -\| y - AFx \|^2 / 2\sigma^2 - \beta \text{TV}(x) - 1_{\left[0,1\right]^d}(x) \right] ,
\end{equation}
with $\beta = 100$. {Above, $1_{[0,1]^d}$ is the convex indicator function on the unit cube, as the pixel values for this experiment are scaled to the range $[0,1]$.}
Following~{(\ref{eq:g_gbar}) and (\ref{eq:BigF})}, we have that $f(x)=\| y - AFx \|^2 / 2\sigma^2$ and $g(x) = {\beta}\text{TV}(x) + 1_{\left[0,1\right]^d}(x)$.
Figure \ref{fig:tomography_x_y}(a) shows the {Shepp-Logan phantom} tomography test image for this experiment and Figure \ref{fig:tomography_x_y}(b) shows the amplitude of the {(noisy)} Fourier coefficients 
{collected in the}
observation {vector} $y$ (in logarithmic scale). {In this figure,} black regions represent unobserved pixels.

\begin{figure}
    \centering
    \subfigure[$\log\pi^{\gamma}(x)$]
    {
        \includegraphics[scale=.25]{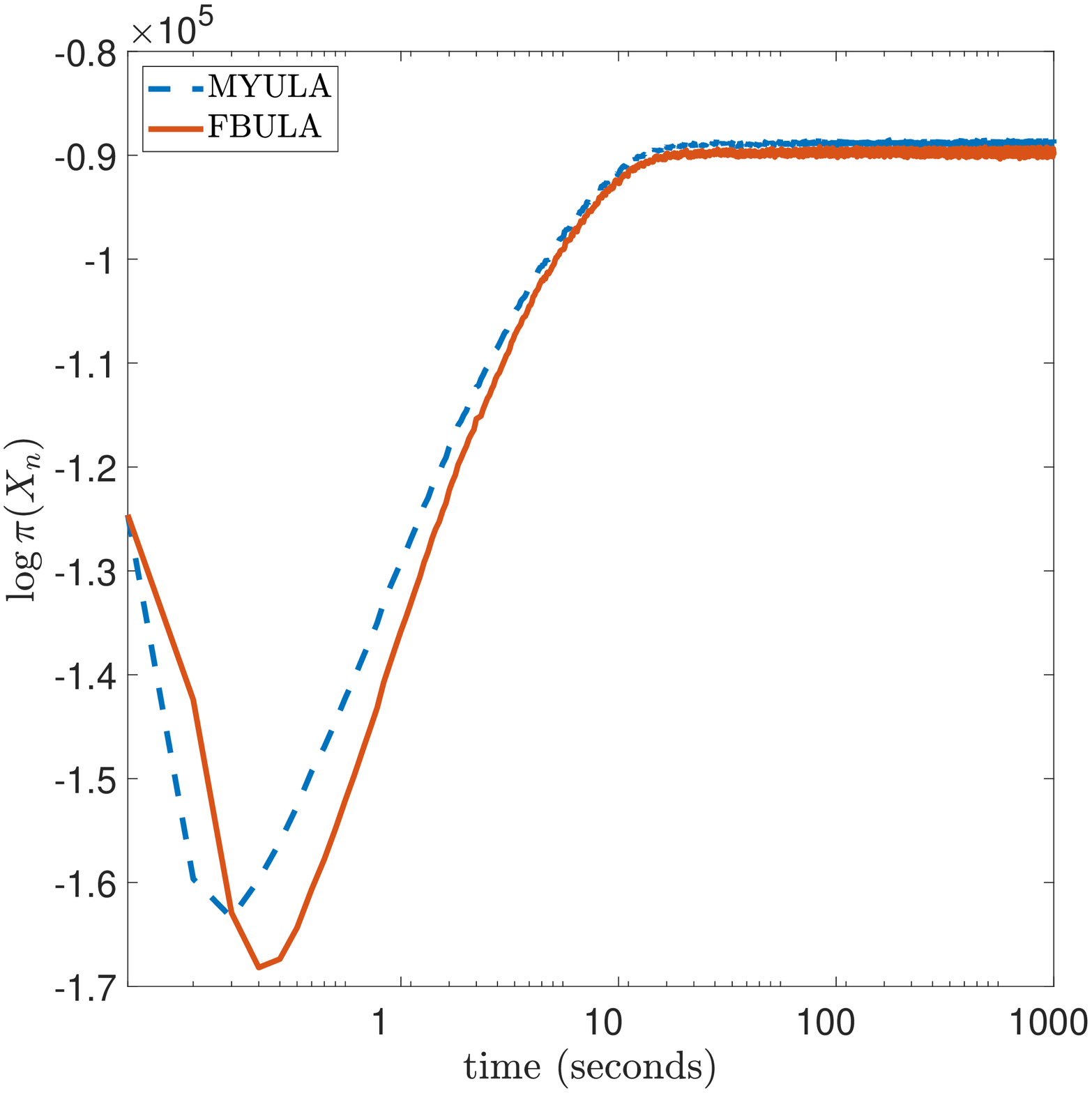}
        \label{subfig:tomography_logPiTrace}
    }
    \subfigure[MSE]
    {
        \includegraphics[scale=.25]{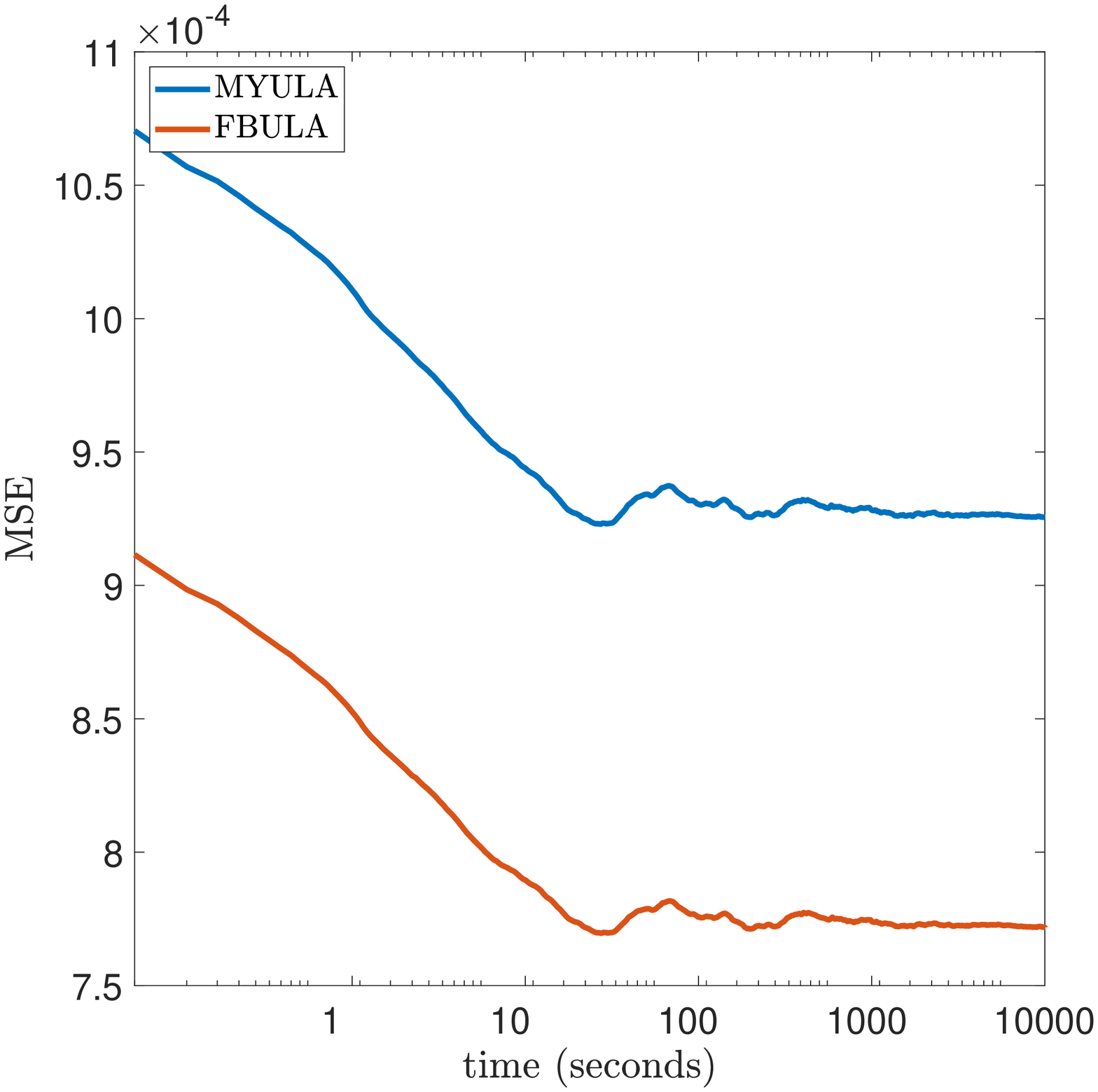}
        \label{subfig:tomography_mse} 
    }\\
    \subfigure[MYULA: posterior mean.]
    {
        \includegraphics[scale=.75]{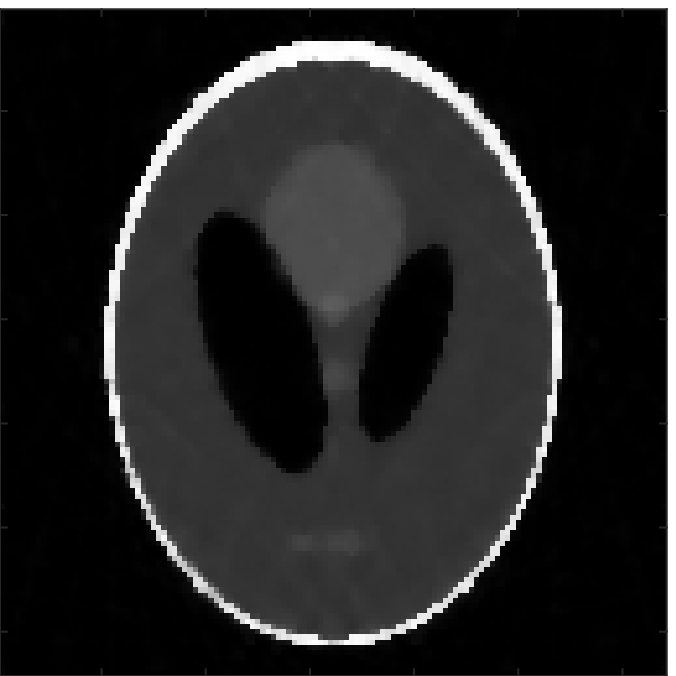}
        \label{subfig:tomography_myula_mean}
    }
    \subfigure[FBULA: posterior mean.]
    {
        \includegraphics[scale=.75]{tomography_post_mean_FBULA}
        \label{subfig:tomography_fbula_mean} 
    }
    \caption{\textit{Tomography} experiment: (a) Convergence to the typical set of the posterior distribution (\ref{eq:tomography_posterior_dist}). (b) Evolution of the MSE in stationarity. Posterior mean of (\ref{eq:tomography_posterior_dist}) as estimated with (c) MYULA and (d) FBULA, { respectively, after} $10^6$ {iterations}.
    }
    \label{fig:tomography_logPiTrace_mse}
\end{figure}

We have set $\gamma = 1/5 L_f$ where $L_f = 1/\sigma^2 = 10^4$ {for both MYULA and FBULA.} Figure~\ref{fig:tomography_logPiTrace_mse}(a) shows the evolution of the values of $\log \pi(x)$ from (\ref{eq:tomography_posterior_dist}) of {both MYULA and FBULA}
with {the} step-size $h=1/(L_f + 1/\gamma) = 1.67 \times 10^{-5}$. We {observe}
that both methods converge at a similar rate. However, Figure~\ref{fig:tomography_logPiTrace_mse}(b) shows the evolution of the mean-squared error (MSE) 
between the {ergodic} mean of the samples 
and the true image~$x$. 
{Here, it}
can be seen  FBULA reaches a better MSE level 
{compared to}
 MYULA.  We have also included  in Figure~\ref{fig:tomography_logPiTrace_mse}(c),(d) the {posterior mean estimated by both MYULA and FBULA.}

\bibliographystyle{unsrt}
\bibliography{references}

\appendix 

\renewcommand{\theequation}{\thesection.\arabic{equation}}
\setcounter{equation}{0}


\section{Proof of Proposition \ref{prop:my}}

To establish the admissibility of the MY envelopes, we verify the conditions in Definition \ref{assump:envAssumption}.
We begin by verifying that Definition \ref{assump:envAssumption}\ref{item:gradLips} is met. 
Recall that the MY envelope is continuously-differentiable \cite[Theorem 2.26]{rockafellar2009variational},
and its gradient is given by
\begin{align}
    \nabla F_\g^{\MY}(x) = \nabla f(x) + \frac{x-P_{\g g}(x)}{\g}, \qquad x\in \R^d,
    \label{eq:MYGrad}
\end{align}
where 
\begin{equation}
    P_{\g g}: x \rightarrow \argmin_{z\in \R^d}\left\{g(z)+\frac{1}{2\gamma}\|x-z\|_2^{2} \right\}
    \label{eq:proxDefn}
\end{equation}
is the proximal operator associated with the function $\g g$. Above, note that the  minimizer is unique and the map $P_{\g g}$ is thus well-defined. 
Recall also that we can use the Moreau decomposition~\cite{beck2017first} to write that 
\begin{equation}
    x = P_{\g g}(x) + P_{\g g^*(\cdot/\g)}(x), \qquad x\in \R^d,
    \label{eq:moreauDec}
\end{equation}
where $g^*$ is the Fenchel conjugate of $g$. Using \eqref{eq:moreauDec}, we rewrite \eqref{eq:MYGrad} as
\begin{equation}
    \nabla F^{\MY}_\g(x)  = \nabla f(x) + \frac{P_{\g g^*(\cdot/\g)}(x)}{\g},
    \label{eq:eqGrad}
\end{equation}
which we will use next to compute the Lipschitz constant of $\nabla F_\g^{\MY}$.
For $x,y\in \R^d$, note that 
\begin{align}
    & \|\nabla F_\g^{\MY}(x) - \nabla F_\g^{\MY}(y)\|_2 \nonumber\\
    & \le \| \nabla f(x) - \nabla f(y) \|_2 + \frac{1}{\g} \| P_{\g g^*(\cdot/\g)}(x) - P_{\g g^*(\cdot/\g)}(y)\|_2
    \qquad \text{(see \eqref{eq:eqGrad})}
    \nonumber\\
    & \le \lambda_2 + \frac{1}{\g},
    \qquad \text{(\eqref{eq:gradLipschitz} and non-expansiveness of the proximal operator)}
\end{align}
where the first line above uses the triangle inequality. 
We next verify that Definition \ref{assump:envAssumption}\ref{item:pntCvg} holds. 
In one direction, it is easy to see that $g_\g\le g$ for every $\g>0$.  In the other direction, we 
fix $x\in \R^d$ and distinguish two cases: 
\begin{enumerate}[leftmargin=*,wide]
    \item When $x\in \CC$, let us fix an arbitrary $\epsilon>0$. Because $g(x)<\infty$ and $\min_z g(z) > -\infty$ by \eqref{eq:propsG}, there exists a sufficiently small~$\g_\epsilon>0$ such that the following holds for every $\g\le \g_\epsilon$:
    \begin{equation}
        \min_{z\in \R^d}\l\{ g(z) + \frac{\epsilon^2}{2\g}\r\} > g(x). 
        \qquad x\in \CC.
        \label{eq:suffSmallGamma}
    \end{equation}
    We now use the above inequality as part of the following argument, 
    \begin{align}
        & \min_{z\in \R^d}\l\{ g(z)+ \frac{1}{2\g}\|x-z\|_2^2: \|x-z\|\ge \epsilon \r\} \nonumber\\
        & \ge \min_z \l\{g(z)+ \frac{\epsilon^2}{2\g}\r\} \nonumber \\
        & > g(x) \qquad \text{(see \eqref{eq:suffSmallGamma})} \nonumber\\
        & \ge g_\g(x) = \min_z \l\{ g(z)+ \frac{1}{2\g}\|x-z\|_2^2  \r\},
        \label{eq:chain}
    \end{align}
which holds when $x\in\CC$ and for every $\g\le \g_\epsilon$. In view of \eqref{eq:proxDefn} and \eqref{eq:chain}, we conclude that $\|x-P_{\g g}(x)\|_2\le \epsilon$, for every $x\in \CC$ and every $\g \le \g_\epsilon$. Since the choice of $\epsilon$ was arbitrary, we arrive at 
\begin{equation}
\lim_{\g\rightarrow0}\|x - P_{\g g}(x)\|_2 = 0,
\qquad x\in \CC.
\label{eq:insideC}
\end{equation}
It immediately follows that 
    \begin{align}
        & \liminf_{\g\rightarrow 0} g_{\g}(x) - g(x) \nonumber\\
        & = \liminf_{\g\rightarrow 0} \min_{z}\l\{g(z)+\frac{1}{2\g}\|x-z\|_2^2 \r\} - g(x) 
        \qquad \text{(see \eqref{eq:MYEnv})}
        \nonumber\\
        & = \liminf_{\g\rightarrow 0} \l\{g(P_{\g g}(x)) + \frac{1}{2\g}\|x - P_{\g g}(x)\|_2^2\r\} - g(x) 
        \qquad  \text{(see \eqref{eq:proxDefn})}
        \nonumber\\
        & \ge \liminf_{\g \rightarrow 0} g(P_{\g g}(x)) - g(x) \nonumber\\
        & \ge \liminf_{z\rightarrow x} g(z) - g(x) 
        \qquad \text{(see \eqref{eq:insideC})}\nonumber\\
        & \ge 0, \qquad x\in \CC,
        \label{eq:limInf}
    \end{align}
    where the last inequality follows from the fact that $g$ is lower semi-continuous, see above Equation~\eqref{eq:propsG}. If we combine \eqref{eq:limInf} with the earlier observation that $g_\g \le g$ for every $\g$, we reach the conclusion that~$\lim_{\g\rightarrow 0} g_\g(x)=g(x)$, provided that $x\in \CC$. Lastly, after recalling the definitions of $F$ and $F_\g^\MY$ in \eqref{eq:BigF} and~\eqref{eq:MYEnv}, respectively, we arrive at    
    \begin{equation}
        \lim_{\g\rightarrow 0} F^\MY_\g(x) = F(x),
        \qquad x\in \CC.
    \end{equation}

\item When $x\notin \CC$, note that $g(x)=\infty$ by \eqref{eq:propsG}. Note also that 
\begin{align}
    g_\g(x) & = \min_z\l\{ g(z)+ \frac{1}{2\g} \|x-z\|_2^2 \r\}
    \qquad \text{(see \eqref{eq:MY})}
    \nonumber\\
    & = \min \l\{g(z) + \frac{1}{2\g} \|x- z\|^2_2: z\in \CC\r\}\qquad \text{(see \eqref{eq:propsG})}
    \nonumber\\
    & \ge -\max_{z\in \CC} g(z) + \frac{1}{2\g} \min_{z\in \CC} \|x-z\|_2^2 
    \nonumber\\
    & =: -\max_{z \in \CC}g(z)+ \frac{1}{2\g} \dist(x,\CC)^2,
    \label{eq:gToDist}
\end{align}
where $\dist(x,\CC)$ is the Euclidean distance from $x$ to the set $\CC$.  The maximum above is finite by \eqref{eq:propsG}. 
Above, by sending $\g$ to zero, we immediately find that $\lim_{\g\rightarrow 0} g_\g(x)=\infty=g(x)$,  provided that  $x\notin \CC$. Recalling the definition of $F$ and $F_\g$ in \eqref{eq:BigF} and \eqref{eq:MYEnv}, we conclude that  
\begin{equation}
\lim_{\g\rightarrow 0} F_\g^\MY(x) = F(x) = \infty,\qquad  x\notin \CC.
\label{eq:outsideC}
\end{equation}
\end{enumerate}
Together, \eqref{eq:insideC} and \eqref{eq:outsideC} imply that $\lim_{\g\rightarrow 0} F_\g(x) = F(x)$ for every $x\in \R^d$. 

Lastly, we now verify  that Definition~\ref{assump:envAssumption}\ref{item:integrable} is satisfied: Because, by Assumption \ref{assump:envAssumption}\ref{item:CAssumption}, $\CC$ is enclosed inside a ball of radius~$R$ centered at the origin, it holds that 
\begin{align}
    \dist(x,\CC) \ge (\|x\|-R)_+, \qquad x\in \R^d,
    \label{eq:lowBndDist}
\end{align}
where $(a)_+:=\max(a,0)$. When  $\|x\|_2$ is sufficiently large, we can simplify \eqref{eq:lowBndDist}. In particular,~\eqref{eq:lowBndDist} immediately implies that 
\begin{equation}
    \dist(x,\CC) \ge \|x\|_2/2, \qquad \|x\|_2 \ge 2R.
    \label{eq:gToDist2}
\end{equation}
If $\|x\|_2\ge 2R$, then we can use the convexity of $f$ by Assumption \ref{assumption:fg}\ref{item:fAssump}, in order  to write that 
\begin{align}
    F_\g^\MY(x) 
    &= f(x)+g_\g(x) \qquad \text{(see \eqref{eq:MYEnv})} \nonumber\\
    & \ge  f(0)+\langle x, \nabla f(0)\rangle + g_\g(x)
    \qquad (\text{convexity of }f)
    \nonumber\\
    & \ge f(0)+\langle x, \nabla f(0) \rangle - \max_{z\in \CC}g(z) + \frac{1}{2\g}\dist(x,\CC)^2
    \qquad \text{(see \eqref{eq:gToDist})}
    \nonumber\\
    & \ge f(0)+\langle x, \nabla f(0) \rangle - \max_{z\in \CC}g(z) + \frac{\|x\|_2^2}{8\g},
    \qquad \text{(see \eqref{eq:gToDist2})}
\end{align}
for $ \|x\|_2 \ge 2R$. 
We now set 
\begin{align}
    F^0(x) := 
    \begin{cases}
    F_\g(x) & \|x\|_2 \le 2R \\
    f(0)+\langle x, \nabla f(0) \rangle - \max_{z\in \CC}g(z) + \frac{\|x\|_2^2}{8\g} & \|x\|_2 \ge 2R.
    \end{cases}
\end{align}
By its construction, note that $F^0$ satisfies  $\int e^{-F^0(x)}\, \der x<\infty$ and $F^0\ge F_\g$,
as required in Definition~\ref{assump:envAssumption}\ref{item:integrable}. The former claim is true because $F^0(x)$ is quadratic for large $\|x\|_2$ and $e^{-F^0}$ thus decays rapidly faraway from the origin.
This completes the proof of Proposition~\ref{prop:my}.



\section{Proof of Proposition \ref{prop:fwbwAdmissible}}

To establish the admissibility of the FB envelopes, we verify the requirements in Definition~\ref{assump:envAssumption}. 
We first verify that Definition \ref{assump:envAssumption}\ref{item:gradLips} holds. To begin, for $x\in \R^d$, let 
\begin{equation}
\ol{x} := x-\g \nabla f(x), \qquad T_\g(x):= P_{\g g}(\ol{x}),\label{eq:shortHands}
\end{equation}
for short. Above, recall that $P_{\g g}$ is the proximal operator associated with the function $\g g$, see~\eqref{eq:gradMY+proximal}. 
We can then apply the Moreau decomposition~\cite{beck2017first} to write that 
\begin{equation}
    \ol{x} = P_{\g g}(\ol{x}) + P_{\g g^*(\cdot/\g)}(\ol{x}),\label{eq:2ndMeuro}
\end{equation}
where $g^*$ is the Fenchel conjugate of $g$. 
For future reference, we record the following observations:
\begin{align}
    \|T_\g(x)\|_2 & = \|P_{\g g}(\ol{x}) \|_2 \le R,
    \qquad \text{(combine Assumption~\ref{assumption:fg}\ref{item:CAssumption} and 
    \eqref{eq:gToDist})}
    \label{eq:normTx}
\end{align}
\begin{align}
        \|x-T_\g(x)\|_2 & \le \|x\|_2 + \|T_\g(x)\|_2 \le \|x\|_2+ R, \qquad \text{(see \eqref{eq:normTx})}
        \label{eq:x-Tx}
\end{align}
\begin{align}
    &\| x -y - (T_\g(x) -T_\g(y))\|_2  \nonumber\\
    & =\| x-y - (P_{\g g}(\ol{x}) - P_{\g g}(\ol{y}))\|_2 
    \qquad \text{(see \eqref{eq:shortHands})}
    \nonumber\\
    & = 
    \| \g (\nabla f(x)-\nabla f(y)) + (P_{\g g^*(\cdot/\g)}(\ol{x})-P_{\g g^*(\cdot/\g)}(\ol{x}))\|_2 
    \quad \text{(see \eqref{eq:2ndMeuro})}
    \nonumber\\
    & \le \g \| \nabla f(x) - \nabla f(y)\|_2
     + \l\| P_{\g g^*(\cdot/\g)}(\ol{x})-P_{\g g^*(\cdot/\g)}(\ol{x})\r\|_2 \nonumber\\
     & \le \g \lambda_2 \| x-y\|+ \| \ol{x}-\ol{y}\|_2 
     \qquad (\text{\eqref{eq:gradLipschitz} and non-expansiveness of }P_{\g g^*(\cdot/\gamma)})
     \nonumber\\
     & \le \g\lambda_2 \|x-y\|_2+ \|x-y \|_2+ \g \|\nabla f(x) - \nabla f(y)\|_2 
    \qquad \text{(see \eqref{eq:shortHands})}
     \nonumber\\
     & \le (1+2\g \lambda_2) \| x-y\|,
     \qquad \text{(see \eqref{eq:gradLipschitz})}
     \label{eq:x-y-Tx+Ty}
\end{align}
\begin{align}
  & 0 \prec (1- \g \lambda_2)I_d \preccurlyeq I_d - \g \nabla^2f(x) \preccurlyeq I_d, \nonumber\\
&
    \| I_d - \g \nabla^2 f(x)\| \le 1. \qquad (\text{\eqref{eq:gradLipschitz} and } \g\in (0,1/\lambda_2))
    \label{eq:specofR}
\end{align}
Above, $I_d\in \R^{d\times d}$ is the identity matrix. 
Recall the expression for $\nabla F_\g^{\FB}$  from Proposition~\ref{prop:propsFB}\ref{item:derFB}. We next establish that $F_\g^\FB$ is smooth. Below, without loss of generality, we can assume that $\|y\|_2\le \lambda_0$: 
\begin{align}
    & \|\nabla F^\FB_\g(x) - \nabla F^\FB_\g(x)\|_2 \nonumber\\
    & = \frac{1}{\g} \| (I_d-\g \nabla^2 f(x)) (x - T_\g(x) ) - (I_d-\g \nabla^2 f(y)) (y - T_\g(y) ) \|_2
    \qquad \text{(Proposition~\ref{prop:propsFB}\ref{item:derFB})}
    \nonumber\\
    & \le  \frac{1}{\g}\l\| (I_d-\g \nabla^2 f(x)) (x-y - (T_\g(x)-T_\g(y)) \r\|_2 
    +  \l\|(\nabla^2 f(x) -\nabla^2 f(y)) (y-T_\g(y))  \r\|_2 \nonumber\\
    & \le \frac{1+2\g\lambda_2}{\g} \|x-y\|_2
    + \lambda_3 (\|y\|_2+R) \|x-y\|_2 
    \qquad \text{(see \eqref{eq:hessianLipschitz} and  \eqref{eq:x-Tx}-\eqref{eq:specofR})}
    \nonumber\\
    & \le \frac{1+2\g\lambda_2}{\g} \|x-y\|_2
    + \lambda_3 (\lambda_0+R) \|x-y\|_2  
    \qquad (\|y\|_2\le \lambda_0 \text{ without loss of generality})
    \nonumber\\
    & = (\g^{-1}+2\lambda_2+\lambda_3(\lambda_0+R))\|x-y\|_2 \nonumber\\
    & =: \lambda_{\g}^{\FB} \|x-y\|_2. 
    \label{eq:gradFBLips}
\end{align}
Indeed, if both $\|x\|_2\ge \lambda_0$ and $\|y\|_2 \ge \lambda_0$, then $f(x)=f(y)=\|\nabla f(x)\|_2=\|\nabla f(y)\|_2=\|\nabla^2f(x)\|=\|\nabla^2f(y)\|=0$ by \eqref{eq:lipsAll}, and \eqref{eq:gradFBLips} still holds. This observation justfies our earlier restriction to the case where $\|y\|_2\le \lambda_0$. 
Next, we show the strong convexity of $F_\g^\FB$ over long distances. Again, without loss of generality, we assume below that $\|y\|_2\le \lambda_0$ and write that 
\begin{align}
    & \langle x-y, \nabla F_\g^\FB(x) - \nabla F_\g^\FB(y) \rangle\nonumber\\
    & = \frac{1}{\g} \l\langle x-y, (I_d-\g \nabla^2 f(x)) (x-T_\g(x))- (I_d-\g \nabla^2 f(y)) (y-T_\g(y)) \r\rangle  
    \nonumber\\
    & = \frac{1}{\g} \l\langle x-y, (I_d-\g \nabla^2f(x)) ( x-y - (T_\g(x)-T_\g(y)) ) \r\rangle\nonumber\\
    & \qquad - \l\langle x-y, (\nabla^2f(x)-\nabla^2 f(y)) (y-T_\g(y)) \r\rangle 
    \nonumber\\
    & \ge  \frac{1}{\g}\langle x-y,(I_d-\g \nabla^2 f(x)) (x-y)\rangle - \frac{\|I_d-\g \nabla^2 f(x)\|}{\g} \|x-y\|_2 (\|T_\g(x)\|_2+\|T_\g(y)\|_2)
    \nonumber\\
    & \qquad - \|x-y\|_2 \cdot \|\nabla^2 f(x) -\nabla^2f(y)\|\cdot  \|y-T_\g(y)\|_2 
    \nonumber\\
    & \ge \frac{1-\g\lambda_2}{\g}\|x-y\|_2^2 \nonumber\\
    & \qquad - \frac{2R}{\g} \|x-y\|_2 -\lambda_3(\lambda_0+R) \|x-y\|_2^2
    \qquad \text{(see \eqref{eq:hessianLipschitz},  \eqref{eq:normTx}, \eqref{eq:x-Tx}, \eqref{eq:specofR})}
    \nonumber\\
    & = \frac{1 }{\g}(1- \g( \lambda_2+ \lambda_3(\lambda_0+R)))\|x-y\|_2^2 - \frac{2R}{\g}\|x-y\|_2,
\end{align}
where the second line uses Proposition \ref{prop:propsFB}\ref{item:derFB}. 
When $\|x-y\|_2$ is sufficiently large, the first term in the last line above dominates the second term. In particular, it holds that
\begin{align}
     & \|x-y\|_2 \ge \rho_\g^\FB:= \frac{2R}{1-2\g( \lambda_2+\lambda_3(\lambda_0+R)))}
     \nonumber\\
     &   \qquad \Longrightarrow \langle x-y, \nabla F_\g^\FB(x) - \nabla F_\g^\FB(y) \rangle \ge  ( \lambda_2+ \lambda_3(\lambda_0+R)) \|x-y\|_2^2.
\end{align}
In words, $F^\FB_\g$ behaves like a strongly convex function over long distances. 
It is not difficult to verify that Definition~\ref{assump:envAssumption}\ref{item:pntCvg} is also valid for the FB envelopes: To that end, one needs to combine Proposition~\ref{prop:my} with the relation between the MY and FB envelopes in Proposition~\ref{prop:propsFB}\ref{item:sandwitch2}.

Lastly, we now verify that the requirement in Definition \ref{assump:envAssumption}\ref{item:integrable} is met for the FB envelopes: Below, suppose that $\g\in (0,1/\lambda_2)$. On the one hand, Proposition~\ref{prop:propsFB}\ref{item:sandwitch2} implies that $$F^\FB_\g \ge F^\MY_{\frac{\g}{1-\g\lambda_2}}.$$ On the other hand, by Proposition~\ref{prop:my}, there exists $F^0\le F_{\frac{\g}{1-\g\lambda_2}}^\MY$ such that $e^{-F^0}$ is integrable. By combining the preceding two observations, we find that the FB envelopes satisfy Definition \ref{assump:envAssumption}\ref{item:integrable} for every $\g\in (0,1/\lambda_2)$. This completes the proof of Proposition~\ref{prop:fwbwAdmissible}.  


\section{Proof of Theorem \ref{prop:distMeasures}}

Let  us define 
\begin{equation}
    \Q(x) := \frac{e^{f(x)-1_\CC(x)}}{\int_z e^{f(z)-1_\CC(z)}\, \der z}, \qquad x\in \R^d,
    \label{eq:Q0}
\end{equation}
where $1_\CC$ is the indicator function on the set $\CC$. Note that the target distribution~$\P$ in~\eqref{eq:pi} coincides with~$\Q$ in the special case where $g=1_\CC$. 
Let $\dist(x,\CC)$ denote the distance from $x$ to the set $\CC$. For~$\g>0$, we also define
\begin{align}
    & 1_{\CC,\g}(x) := \frac{1}{2\g}\dist(x,\CC)^2 = \frac{1}{2\g}\min_{z\in \CC} \|x-z\|_2^2 \nonumber\\
    & = \min_{z\in \R^d} \l\{ 1_\CC(x) + \frac{1}{2\g}\|x-z\|_2^2 \r\}, \qquad x\in \R^d, 
    \label{eq:MYenvOfInd}
\end{align}
to be the MY envelope of the indicator function on the set $\CC$. We  denote the corresponding probability distribution by 
\begin{equation}
\Q_\g(x) := \frac{e^{f(x)-1_{\CC,\g}(x)}}{\int_z e^{f(z)-1_{\CC,\g}(z)} \, \der z}, \qquad x\in \R^d.
\label{eq:Qgamma}
\end{equation}
When $\g$ is sufficiently small, we may intuitively regard $\Q_\g$ as a proxy for $\Q$. We also recall our earlier notation for the convenience of the reader: 
\begin{align}
    \P(x) = \frac{e^{-F(x)}}{\int_z e^{-F(z)}\,\der z},
    \qquad 
    \P^\FB_\g(x) = \frac{e^{-F^\FB_\g(x)}}{\int_z e^{-F^\FB_\g(z)}\,\der z}, 
    \qquad x\in \R^d.
\end{align}
Above, the functions $F$ and $F^\FB_\g$ were defined in \eqref{eq:BigF} and \eqref{eq:fbEnv}, respectively. Our objective is to control the distance $\W_1(\P^\FB_\g,\P)$. To begin, we  use the triangle inequality to write that 
\begin{equation}
    \W_1(\P^\FB_\g,\P) \le \W_1(\P^\FB_\g,\Q_\g) + \W_1(\Q_\g,\Q) + \W_1(\Q,\P). 
    \label{eq:triangle3}
\end{equation}
The following three lemmas each bounds one of the terms on the right-hand side above. 

\begin{lem}\label{lem:1stLeg}It holds that 
\begin{equation}
    \W_1(\P^\FB_\g,\Q_\g)
     \le \Cr{dist1} + \Cr{dist2} \frac{I_2(\g/(1-\g\lambda_2)))}{\vol(\CC)+I_1(\g)},\label{eq:1stlegEq}
\end{equation}
where 
\begin{align}
    & \Cr{dist1}:=  \frac{e^{2\max_{x\in \CC} g(x)-\g\lambda_2 \min_z f(z)} \int_{\CC} \|x\|_2 e^{-(1-\g\lambda_2)f(x)} \, \der x }{\int e^{-f(x)-\frac{\dist(x,\CC)^2}{2\g}}\,\der x}\nonumber\\
    & \qquad - \frac{e^{\g\lambda_2 \min_z f(z)-2\max_{x\in \CC} g(x)} \int_\CC \|x\|_2 e^{-f(x)}\, \der x}{\int e^{-(1-\g\lambda_2)f(x) - \frac{\dist(x,\CC)^2}{2\g/(1-\g\lambda_2)}}\, \der x},\nonumber\\
    & \Cr{dist2}:= e^{\max_x f(x) -  \min_x f(x) +2\max_{x\in \CC}g(x)}, \nonumber\\
    & I_1(\g) :=  \sum_{i=0}^{d-1} \vol_{i}(\CC) \cdot (2\pi \g)^{\frac{d-i}{2}}, \nonumber\\
    & I_2(\g):= \sum_{i=0}^{d-1} \vol_i(\CC)\cdot  (2\pi\g)^{\frac{d-i}{2}} \l( \sqrt{\g(d-i+3)}+R \r) =: I_2(\g),
\end{align}
where $\vol_i(\CC)$ is the $i$-th intrinsic volume of $\CC$~\cite{klain1997introduction}. In particular, the $d$-th volume of $\CC$ coincides with the standard volume of $\CC$, i.e.,  $\vol_d(\CC)=\vol(\CC)$. 

\end{lem}

\begin{lem}\label{lem:2ndLeg}It holds that
\begin{equation}
    \W_1(\Q_\g,\Q) \le 
    \Cr{dist3}\frac{R I_1(\g)+ I_2(\g)}{\vol(\CC)+I_1(\g)},
    \label{eq:2ndLegEq}
\end{equation}
where $\Cr{dist3}:=e^{\max_x f(x)-\min_x f(x)}$.
\end{lem}
\begin{lem}\label{lem:3rdLeg}It holds that 
\begin{equation}
    \W_1(\Q,\P) \le \Cr{dist4} R,
    \label{eq:3rdLegEq}
\end{equation}
where 
\begin{equation}
    \Cr{dist4}:= e^{\max_{x\in \CC}g(x)-\min_{x\in \CC}g(x) }-e^{\min_{x\in \CC}g(x)-\max_{x\in \CC} g(x)}.
\end{equation}

\end{lem}
To keep the notation light, above we have suppressed the dependence of $\Cr{dist1}$ to $\Cr{dist4}$ on $\CC,f,g$.
With Lemmas~\ref{lem:1stLeg}-\ref{lem:3rdLeg} at hand, we revisit \eqref{eq:triangle3} and write that 
\begin{align}
    \W_1(\P_\g^\FB,\P) \le 
   \Cr{dist1}+\frac{
    \Cr{dist3} R I_1(\g)+\Cr{dist3} I_2(\g) + \Cr{dist2} I_2(\g/(1-\g\lambda_2))
   }{\vol(\CC)+I_1(\g)}+   \Cr{dist4} R.
\end{align}
This completes the proof of Proposition \ref{prop:distMeasures}.

\section{Proof of Lemma \ref{lem:1stLeg}}

In order to upper bound the distance $\W_1(\Q_\g,\P_\g^\FB)$, we will use the following simple result which bounds the $\W_1$ distance under small perturbations. 
\begin{lem}\label{lem:h1h2lemma}
For constants $\a\ge \b$ and $\b'\le 1$, consider two functions $h_1:\R^d\rightarrow\R$ and $h_2:\R^d\rightarrow\R$ that are related as 
\begin{align}
    & \b' \cdot h_1(x)+\b\le h_2(x) \le h_1(x)+\a, \qquad x\in \R^d,
    \label{eq:h1h2}
\end{align}
Then it holds that 
\begin{align}
    & \W_1\l( \frac{e^{-h_1}}{\int e^{-h_1}}, \frac{e^{-h_2}}{\int e^{-h_2}} \r) \nonumber\\
    & \le 
    \frac{e^{\a-\b} \int \|x\|_2 e^{-\b' h_1(x)}\,\der x}{\int e^{-h_1(x)}\,\der x}
    - \frac{e^{\b-\a}\int \|x\|_2 e^{-h_1(x)}\, \der x}{\int e^{-\b' h_1(x)}\, \der x}.
    \label{eq:h1h2result}
\end{align}
\end{lem}
As a sanity check, note that the right-hand side of \eqref{eq:h1h2result} is always nonnegative because, by assumption, $\a\ge \b$ and $\b'\le 1$. Moreover, if we set $\a=\b=0$ and $\b'=1$, the right-hand side of \eqref{eq:h1h2result} reduces to zero, as expected.

Let us now recall from \eqref{eq:Qgamma} and \eqref{eq:pi} that $\Q_\g \propto e^{-f-1_{\CC,\g}}$ and $\P_\g^\FB\propto e^{-F^\FB_\g}$, respectively.  Our plan is to invoke Lemma~\ref{lem:h1h2lemma} with the choice of $h_1=f+1_{\CC,\g}$ and $h_2=F^\FB_\g$. In turn, this plan necessitates that we verify~\eqref{eq:h1h2} for the choice of $h_1=f+1_{\CC,\g}$ and $h_2=F^\FB_\g$. 
We begin by relating the two functions as 
\begin{align}
    F^\FB_\g(x) & \le F^\MY_\g(x) 
    \qquad \text{(Proposition \ref{prop:propsFB}\ref{item:sandwitch2})}\nonumber\\
    & = f(x) + g_\g(x) \qquad \text{(see \eqref{eq:MYEnv})} \nonumber\\
    & = f(x) + \min_{z\in \CC}\l\{g(z) + \frac{1}{2\g}\|x-z\|_2^2\r\}
    \qquad \text{(see \eqref{eq:MY})}\nonumber\\
    & \le \max_{z\in \CC}g(z) + f(x)+ \frac{1}{2\g} \min_{z\in \CC} \|x-z\|_2^2\nonumber\\
    & =\max_{z\in \CC}g(z)+f(x) + 1_{\CC,\g}(x)\qquad \text{(see \eqref{eq:MYenvOfInd})}.
\end{align}
Above, note that $\max_{z\in \CC} g(z)$ is finite by the construction of $g$, see~\eqref{eq:propsG}. In the other direction, we similarly write that 
\begin{align}
    & F^\FB_\g(x) \nonumber\\
    & \ge F^\MY_{\frac{\g}{1-\g\lambda_2}}(x) 
    \qquad \text{(Proposition \ref{prop:propsFB}\ref{item:sandwitch2})}\nonumber\\
    & = f(x)+g_{\frac{\g}{1-\g\lambda_2}}(x) \nonumber\\
    & = f(x)+\min_z\l\{ g(z)+ \frac{1-\g\lambda_2}{2\g}\|x-z\|_2^2 \r\} \nonumber\\
    & \ge  - \max_{z\in \CC}g(z) + f(x)+ \frac{1-\g\lambda_2}{2\g}\min_{z\in \CC} \|x-z\|_2^2 \nonumber\\
    & =  - \max_{z\in \CC}g(z)
    +\g\lambda_2 f(x)
    +(1-\g\lambda_2)f(x) + (1-\g\lambda_2) 1_{\CC,\g}(x)
    \qquad \text{(see \eqref{eq:MYenvOfInd})} \nonumber\\
    & \ge - \max_{z\in \CC}g(z)
    +\g\lambda_2 \min_z f(z)
    +(1-\g\lambda_2)f(x) + (1-\g\lambda_2) 1_{\CC,\g}(x).
\end{align}
Above, $\min_z f(z)$ is  finite by the construction of $f$, see \eqref{eq:lipsAll}.  
To summarize, for our choice of $h_1=f+1_{\CC,\g}$ and $h_2=F^\FB_\g$, \eqref{eq:h1h2} is satisfied with
\begin{align}
   & \b' (f(x)+1_{\CC,\g}(x))+\b \le F^\FB_\g(x) \le f(x)+1_{\CC,\g}(x)+\a, \qquad x\in \R^d,
    \label{eq:summary1}\\
    & \a = \max_{z\in \CC}g(z),\quad \b = \g\lambda_2 \min_z f(z) -\max_{z\in \CC}g(z),\quad \b'=1-\g\lambda_2.\nonumber
\end{align}
With \eqref{eq:summary1} at hand, we can invoke Lemma \ref{lem:h1h2lemma} to find that 
\begin{align}
     \W_1( \Q_\g,\P_\g^\FB )
    & = \W_1\l( \frac{e^{-f-1_{\CC,\g}}}{\int e^{-f-1_{\CC,\g}}},
    \frac{e^{-F_{\g}^\FB}}{\int e^{-F_{\g}^\FB}}
    \r) \nonumber\\
    & \le  
    \frac{e^{\a-\b} \int \|x\|_2 e^{-\b'f(x)-\b' 1_{\CC,\g}(x)} \, \der x}{\int e^{-f(x)-1_{\CC,\g}(x)}\, \der x} - 
    \frac{ e^{\b-\a} \int \|x\|_2 e^{-f(x)-1_{\CC,\g}(x)} \,\der x}{ \int  e^{-\b'f(x)-\b' 1_{\CC,\g}(x)} \,\der x}
    \nonumber\\
    & = \frac{e^{2\max_{x\in \CC}g(x)-\g\lambda_2 \min_z f(z)} \int \|x\|_2 e^{-(1-\g\lambda_2) f(x)-\frac{\dist(x,\CC)^2}{2\g/(1-\g\lambda_2)} }\,\der x}{\int e^{-f(x)-\frac{\dist(x,\CC)^2}{2\g}}\,\der x}\nonumber\\
    &  -
    \frac{e^{\g\lambda_2 \min_z f(z)-2\max_{x\in \CC}g(x)} \int \|x\|_2 e^{-f(x)-\frac{\dist(x,\CC)^2}{2\g}}\,\der x}{\int e^{-(1-\g\lambda_2) f(x)-\frac{\dist(x,\CC)^2}{2\g/(1-\g\lambda_2)}}\,\der x},
    \label{eq:beforeLemInt}
\end{align}
where the second line uses Lemma \ref{lem:h1h2lemma} and the last line uses \eqref{eq:MYenvOfInd} and  \eqref{eq:summary1}.
To complete the proof, we  need the following result, which is proved in \cite{brosse2017sampling,kampf2009weighted} using the well-known Steiner's formula from integral geometry~\cite{klain1997introduction}. For completeness, a self-contained  proof is given later. 
\begin{lem}\label{lem:integral}
It holds that 
\begin{subequations}
\begin{align}
    & \int_{\CC^c} e^{-\frac{\dist(x,\CC)^2}{2\g} }\, \der x = \sum_{i=0}^{d-1} \vol_{i}(\CC) \cdot (2\pi \g)^{\frac{d-i}{2}}=: I_1(\g),
    \label{eq:int0}\\
    & \int_{\R^d} e^{-\frac{\dist(x,\CC)^2}{2\g} }\, \der x = \sum_{i=0}^d \vol_{i}(\CC) \cdot (2\pi \g)^{\frac{d-i}{2}}= \vol(\CC)+ I_1(\g),
    \label{eq:int1}\\
    & \int_{\CC^c}\|x\|_2 e^{-\frac{\dist(x,\CC)^2}{2\g} }\, \der x \nonumber\\
     & \le \sum_{i=0}^{d-1} \vol_i(\CC)\cdot  (2\pi\g)^{\frac{d-i}{2}} \l( \sqrt{\g(d-i+3)}+R \r) 
    =: I_2(\g),
    \label{eq:int2}
\end{align}\label{eq:integrals}\end{subequations}
where $\CC^c=\R^d\backslash \CC$ is the complement of the set $\CC$ and $\vol_i(\CC)$ is the $i$-th intrinsic volume of $\CC$~\cite{klain1997introduction}. In particular, the $d$-th volume of $\CC$ coincides with the standard volume of $\CC$, i.e.,  $\vol_d(\CC)=\vol(\CC)$. 
\end{lem}
In particular, note that $I_1(\g)$ and $I_2(\g)$ both vanish in the limit of $\g\rightarrow0$, as expected. 
To control the first fraction in the last line of \eqref{eq:beforeLemInt}, we use Lemma \ref{lem:integral} to  write that 
\begin{align}
    & \frac{\int \|x\|_2 e^{-(1-\g\lambda_2) f(x)-\frac{\dist(x,\CC)^2}{2\g/(1-\g\lambda_2)} }\,\der x}{\int e^{-f(x)-\frac{\dist(x,\CC)^2}{2\g}}\,\der x} \nonumber\\
    & = \frac{\int_{\CC} \|x\|_2 e^{-(1-\g\lambda_2)f(x)} \,\der x}{\int e^{-f(x)-\frac{\dist(x,\CC)^2}{2\g}}\, \der x }+ \frac{\int_{\CC^c} \|x\|_2 e^{-(1-\g\lambda_2)f(x) - \frac{\dist(x,\CC)^2}{2\g/(1-\g\lambda_2)}}\, \der x}{\int e^{-f(x)-\frac{\dist(x,\CC)^2}{2\g}}\,\der x} \nonumber\\
    & \le \frac{\int_{\CC} \|x\|_2 e^{-(1-\g\lambda_2)f(x)} \,\der x}{\int e^{-f(x)-\frac{\dist(x,\CC)^2}{2\g}}\, \der x } + \frac{e^{-(1-\g\lambda_2)\min_x f(z)}}{e^{-\max_x f(z)}}\cdot \frac{I_2(\g/(1-\g\lambda_2))}{\vol(\CC)+I_1(\g)},
    \label{eq:1stFrac}
\end{align}
where the last line uses Lemma \ref{lem:integral}. 
For the second fraction in the last line of~\eqref{eq:beforeLemInt}, we similarly use Lemma \ref{lem:integral} to write that 
\begin{align}
    & \frac{ \int \|x\|_2 e^{-f(x)-\frac{\dist(x,\CC)^2}{2\g}}\,\der x}{\int e^{-(1-\g\lambda_2) f(x)-\frac{\dist(x,\CC)^2}{2\g/(1-\g\lambda_2)}}\,\der x}  \nonumber\\
    & =  \frac{ \int_{\CC} \|x\|_2 e^{-f(x)}\,\der x}{\int e^{-(1-\g\lambda_2) f(x)-\frac{\dist(x,\CC)^2}{2\g/(1-\g\lambda_2)}}\,\der x}   + \frac{ \int_{\CC^c} \|x\|_2 e^{-f(x)-\frac{\dist(x,\CC)^2}{2\g}}\,\der x}{\int e^{-(1-\g\lambda_2) f(x)-\frac{\dist(x,\CC)^2}{2\g/(1-\g\lambda_2)}}\,\der x}  \nonumber\\
    & \ge \frac{ \int_{\CC} \|x\|_2 e^{-f(x)}\,\der x}{\int e^{-(1-\g\lambda_2) f(x)-\frac{\dist(x,\CC)^2}{2\g/(1-\g\lambda_2)}}\,\der x} .
    \label{eq:2ndFrac}
\end{align}
We can now use \eqref{eq:1stFrac} and \eqref{eq:2ndFrac} to upper bound the $\W_1$ distance in \eqref{eq:beforeLemInt} as 
\begin{align}
    \W_1(\Q_\g,\P^\FB_\g) & 
    \le \frac{e^{2\max_{x\in \CC} g(x)-\g\lambda_2 \min_z f(z)} \int_{\CC} \|x\|_2 e^{-(1-\g\lambda_2)f(x)} \, \der x }{\int e^{-f(x)-\frac{\dist(x,\CC)^2}{2\g}}\,\der x}\nonumber\\
    & - \frac{e^{\g\lambda_2 \min_z f(z)-2\max_{x\in \CC} g(x)} \int_\CC \|x\|_2 e^{-f(x)}\, \der x}{\int e^{-(1-\g\lambda_2)f(x) - \frac{\dist(x,\CC)^2}{2\g/(1-\g\lambda_2)}}\, \der x} \nonumber\\
    & + e^{\max_x f(x) - \min_x f(x) +2\max_{x\in \CC}g(x)} \frac{I_2(\g/(1-\g\lambda_2)))}{\vol(\CC)+I_1(\g)},
\end{align}
which completes the proof of Lemma \ref{lem:1stLeg}.

\section{Proof of Lemma \ref{lem:h1h2lemma} }

Let us first recall Theorem 6.15 in \cite{villani2009optimal}, which can be used to link the $\W_1$ and TV distances of two rapidly decaying distributions.  
\begin{thm}\label{thm:WtoTV} Any pair $(\Rfrak_1,\Rfrak_2)$ of  probability distributions  satisfies
\begin{equation}
    \W_1(\Rfrak_1,\Rfrak_2) \le \int \|x\|_2\cdot  |\Rfrak_1(x)-\Rfrak_2(x)|\, \der x.
\end{equation}
\end{thm}
We  apply Theorem \ref{thm:WtoTV} to bound the $\W_1$  distance of interest as 
\begin{align}
\W_1\l( \frac{e^{-h_1}}{\int e^{-h_1}}, \frac{e^{-h_2}}{\int e^{-h_2}} \r) \le  \int_{\R^d} \|x\|_2 \l|\frac{e^{-h_1(x)}}{\int e^{-h_1(z)} \,\der z} - \frac{e^{-h_2(x)}}{\int e^{-h_2(z)} \,\der z} \r|\, \der x.
    \label{eq:TVwritten}
\end{align}
To bound the integral on the right-hand side above, we first use the assumed relationship between $h_1$ and~$h_2$ to write that 
\begin{align}
    \frac{e^{-h_1(x)}}{\int e^{-h_1(z)} \,\der z} - \frac{e^{-h_2(x)}}{\int e^{-h_2(z)} \,\der z} & \le \frac{e^{-h_1(x)}}{\int e^{-h_1(z)}\,\der z} - \frac{e^{-h_1(x)-\a}}{ \int e^{-\b' h_1(z)-\b}\,\der z}
    \qquad \text{(see \eqref{eq:h1h2})} \nonumber\\
    & = \frac{e^{-h_1(x)}}{\int e^{-h_1(z)}\, \der z} - \frac{e^{-h_1(x)+\b-\a}}{\int e^{-\b' h_1(z)}\, \der z} \ge 0.\label{eq:abs1}
\end{align}
The last inequality above holds because, by assumption,  $\a\ge \b$ and $\b'\le 1$. In the other direction, we can again use the relation between $h_1$ and $h_2$ to write that 
\begin{align}
     \frac{e^{-h_1(x)}}{\int e^{-h_1(z)} \,\der z} - \frac{e^{-h_2(x)}}{\int e^{-h_2(z)} \,\der z} & \ge \frac{e^{-h_1(x)}}{\int e^{-h_1(z)}\, \der z} - \frac{e^{-\b'h_1(x)-\b}}{\int e^{-h_1(z)-\a}\, \der z}
     \qquad \text{(see \eqref{eq:h1h2})} \nonumber\\
     & = \frac{e^{-h_1(x)}}{\int e^{-h_1(z)}\,\der z}-\frac{e^{-\b' h_1(x)+\a-\b} }{\int e^{-h_1(z)}\, \der z} \le 0,\label{eq:abs2}
\end{align}
where the last inequality above holds because, by assumption,  $\a\ge \b$ and $\b'\le 1$. By combining \eqref{eq:abs1} and~\eqref{eq:abs2}, we find that 
\begin{align}
   & \|x\|_2\cdot  \l| \frac{e^{-h_1(x)}}{\int e^{-h_1(z)} \,\der z} - \frac{e^{-h_2(x)}}{\int e^{-h_2(z)} \,\der z}\r| \nonumber\\
   & \le \|x\|_2\max\l(
    \frac{e^{-\b' h_1(x)+\a-\b}}{\int e^{-h_1(z)}\, \der z}- \frac{e^{-h_1(x)}}{\int e^{-h_1(z)}\, \der z}
    , \frac{e^{-h_1(x)}}{\int e^{-h_1(z)}\, \der z}- \frac{e^{-h_1(x)+\b-\a}}{\int e^{-\b' h_1(z)}\, \der z}\r)
    \nonumber\\
    & \le \|x\|_2 \l( \frac{e^{-\b' h_1(x)+\a-\b}}{\int e^{-h_1(z)}\,\der z}
    - \frac{ e^{-h_1(x)+\b-\a}}{\int e^{-\b' h_1(z)}\, \der z}
    \r),
    \label{eq:W1fixedX}
\end{align}
where the second line uses see \eqref{eq:abs1} and \eqref{eq:abs2}, and the 
last line above uses the inequality $\max(a,b)\le a+b$ for nonnegative scalars $a,b$.
With \eqref{eq:W1fixedX} at hand, we revisit \eqref{eq:TVwritten} and conclude that 
\begin{align}
    \W_1\l( \frac{e^{-h_1}}{\int e^{-h_1}}, \frac{e^{-h_2}}{\int e^{-h_2}} \r) \le 
   \frac{ e^{\a-\b}\int \|x\|_2 e^{-\b' h_1(x)}\, \der x}{\int e^{-h_1(z)}\, \der z} -\frac{ e^{\b-\a} \int \|x\|_2 e^{-h_1(x)}\, \der x  }{\int e^{-\b' h_1(z)}\, \der z},
\end{align}
which completes the proof of Lemma \ref{lem:h1h2lemma}.

\section{Proof of Lemma \ref{lem:integral}}

Let $\iota(A)=1$ if the claim $A$ is true and $\iota(A)=0$ otherwise. Using the fact that 
\begin{equation}
    e^{-\frac{\dist(x,\CC)^2}{2\g}} = \frac{1}{\g}\int_{\dist(x,\CC)}^\infty  te^{-\frac{t^2}{2\g}} \, \der t,
    \label{eq:trivial}
\end{equation}
we can rewrite the left-hand side of \eqref{eq:int1} as 
\begin{align}
    \int_{\R^d} e^{-\frac{\dist(x,\CC)^2}{2\g} }\, \der x & = \frac{1}{\g} \int_{\R^d} \int_{\dist(x,\CC)}^\infty  te^{-\frac{t^2}{2\g}} \, \der t \der x \nonumber\\
    & = \frac{1}{\g}\int_{\R^d} \int_0^\infty \iota(\dist(x,\CC)\le t) t e^{-\frac{t^2}{2\g}}\, \der t\der x \nonumber\\
    & = \frac{1}{\g} \int_0^{\infty} \l(\int_{\R^d} \iota(\dist(x,\CC)\le t) \, \der x \r)\cdot   te^{-\frac{t^2}{2\g}} \, \der t  \nonumber\\
    & = \frac{1}{\g} \int_0^\infty \vol(\{x: \dist(x,\CC)\le t\})\cdot  te^{-\frac{t^2}{2\g}} \der t
\end{align}
Note that $\{x:\dist(x,\CC)\le t\}$ is the tube of radius $t$ around the set $\CC$.  We can represent this set more compactly as $\K+\B(0,t)$, where the addition is in the Minkowski sense. With this in mind, we rewrite the last line above as
\begin{align}
   \int_{\R^d} e^{-\frac{\dist(x,\CC)^2}{2\g} }\, \der x
    & = \frac{1}{\g} \int_0^\infty \vol(\CC + \B(0,t)) \cdot te^{-\frac{t^2}{2\g}} \der t.
    \label{eq:beforeSteiner}
\end{align}
We can now use the Steiner's formula~\cite{klain1997introduction} to express the volume above as
\begin{align}
    \vol(\CC+\B(0,t)) = \sum_{i=0}^d t^i \kappa_i \vol_{d-i}(\CC),
    \label{eq:steiner}
\end{align}
where $\kappa_i$ is the volume of the unit ball in $\R^i$ and $\vol_i(\CC)$ is the $i$-th intrinsic volume of $\CC$~\cite{klain1997introduction}. In particular, the $d$-th intrinsic volume of $\CC$ coincides with its the standard volume, i.e., $\vol_d(\CC)=\vol(\CC)$. Substituting the above identity back into \eqref{eq:beforeSteiner}, we find that
\begin{align}
    & \int_{\R^d} e^{-\frac{\dist(x,\CC)^2}{2\g}}\, \der x \nonumber\\
    & = \frac{1}{\g} \sum_{i=0}^d \kappa_i \vol_{d-i}(\CC) \int_0^\infty t^{i+1}e^{-\frac{t^2}{2\g}} \, \der t 
    \qquad \text{(see \eqref{eq:beforeSteiner} and \eqref{eq:steiner})}
    \nonumber\\
    & = \frac{1}{\g} \sum_{i=0}^d \frac{\pi^\frac{i}{2}}{\Gamma (1+i/2)} \vol_{d-i}(\CC) \cdot \g(2\g)^{\frac{i}{2}}\Gamma(1+i/2) 
    \quad (\kappa_i = \pi^{\frac{i}{2}}/\Gamma(1+i/2) )
    \nonumber\\
    & = \sum_{i=0}^d \vol_{d-i}(\CC) (2\pi\g)^{\frac{i}{2}} \nonumber\\
    & = \sum_{i=0}^d \vol_{i}(\CC) (2\pi \g)^{\frac{d-i}{2}} \nonumber\\
    & = \vol_d(\CC) + \sum_{i=0}^{d-1} \vol_i(\CC) (2\pi\gamma)^{\frac{d-i}{2}} \nonumber\\
    & = \vol(\CC) + \sum_{i=0}^{d-1} \vol_i(\CC) (2\pi\gamma)^{\frac{d-i}{2}}
    \qquad (\vol_d(\CC) = \vol(\CC)) \nonumber\\
    & =: \vol(\CC) +  I_1(\gamma),
\end{align}
which establishes \eqref{eq:int1}. Above, in the second line we used the identity
\begin{equation}
    \int_0^\infty t^i e^{-\frac{t^2}{2\g}} \,\der t = \g(2\g)^{\frac{i-1}{2}} \Gamma\l( \frac{i+1}{2}\r). 
    \label{eq:momentGauss}
\end{equation}
To prove \eqref{eq:int0}, we write that 
\begin{align}
    \int_{\CC^c} e^{-\frac{\dist(x,\CC)^2}{2\g}}\, \der x & = \int_{\R^d} e^{-\frac{\dist(x,\CC)^2}{2\g}} \,\der x- \int_{\CC} e^{-\frac{\dist(x,\CC)^2}{2\g}}\, \der x \nonumber\\
    & = \vol(\CC) + I_1(\gamma) - \int_{\CC} \der x  \nonumber\\
    & = I_1(\g).
\end{align}
To prove \eqref{eq:int2}, we  use the fact that  $\CC\subset \B(0,R)$ by Assumption~\ref{assumption:fg}\ref{item:CAssumption}, which implies that 
\begin{equation}
\dist(x,\CC)\ge \|x\|_2-R.
\label{eq:lwrBndDist}
\end{equation}
In turn, we use \eqref{eq:lwrBndDist} to write that 
\begin{align}
   & \int_{\CC^c} \|x\|_2 e^{-\frac{\dist(x,\CC)^2}{2\g}} \, \der x\nonumber\\
    & \le \int_{\dist(x,\CC)>0} (\dist(x,\CC)+R) e^{-\frac{\dist(x,\K)^2}{2\g}} \, \der x
    \qquad \text{(see \eqref{eq:lwrBndDist})} \nonumber\\
    & = \int_{\dist(x,\CC)>0} (\dist(x,\CC)+R) \l( \int_{\dist(x,\CC)}^\infty \frac{t}{\g}e^{-\frac{t^2}{2\g}}\, \der t \r) \,\der x 
    \qquad \text{(see \eqref{eq:trivial})}
    \nonumber\\
    & \le  \frac{1}{\g} \int_{\R^d}\int_0^\infty \iota(0< \dist(x,\CC) \le t)  (\dist(x,\CC)+R) te^{-\frac{t^2}{2\g}} \, \der t\der x \nonumber\\
    & \le \frac{1}{\g} \int_{\R^d}\int_0^\infty 
    \iota(0< \dist(x,\CC) \le t)
    (t+R)te^{-\frac{t^2}{2\g}}\, \der t\der x \nonumber\\
    & = \frac{1}{\g} \int_0^\infty \l( \int_{\R^d}  \iota(0< \dist(x,\CC) \le t) \,\der x \r)\cdot (t+R)te^{-\frac{t^2}{2\g}} \, \der t \nonumber\\
    & = \frac{1}{\g} \int_0^\infty (\vol( \CC+\B(0,t) )-\vol(\CC)) \cdot (t+R)te^{-\frac{t^2}{2\g}}\, \der t.
    \label{eq:preSteiner2}
\end{align}
We can again use the Steiner's formula to calculate the volume of the tube of radius $t$ around $\CC$, which appears in the last line above. Substituting from \eqref{eq:steiner} into the last line above, we find that 
\begin{align}
    & \int_{\CC^c} \|x\|_2e^{-\frac{\dist(x,\CC)^2}{2\g}} \, \der x \nonumber\\
    & \le  \frac{1}{\g}\int_0^\infty \sum_{i=1}^{d} t^i \kappa_i \vol_{d-i}(\CC) \cdot  (t+R)te^{-\frac{t^2}{2\g}}\, \der t
    \quad (\text{\eqref{eq:steiner}, \eqref{eq:preSteiner2}, }\vol_d(\CC)=\vol(\CC)) \nonumber\\
    & = \frac{1}{\g} \sum_{i=1}^{d}  \kappa_i \vol_{d-i}(\CC) \int_0^\infty (t^{i+2}+Rt^{i+1})e^{-\frac{t^2}{2\g}}\,\der t \nonumber\\
    & = \sum_{i=1}^{d}\kappa_i \vol_{d-i}(\CC) \l( (2\g)^{\frac{i+1}{2}}\Gamma\l(\frac{i+3}{2} \r) + R (2\g)^{\frac{i}{2}} \Gamma\l( \frac{i}{2}+1\r) \r)
    \qquad \text{(see \eqref{eq:momentGauss})} \nonumber\\
    & = \sum_{i=1}^{d} \frac{\pi^{\frac{i}{2}}}{\Gamma(1+i/2)} \vol_{d-i}(\CC) \l( (2\g)^{\frac{i+1}{2}}\Gamma\l(\frac{i+3}{2} \r) + R (2\g)^{\frac{i}{2}} \Gamma\l( \frac{i}{2}+1\r) \r) 
    \quad (\kappa_i = \pi^{\frac{i}{2}}/\Gamma(1+i/2) )
    \nonumber\\
    & = \sum_{i=1}^{d}  \vol_{d-i}(\CC) (2\pi\g)^{\frac{i}{2}} \l(   \sqrt{2\g}\cdot \frac{\Gamma\l( \frac{i+3}{2}\r)}{\Gamma\l( \frac{i+2}{2}\r)} + R  \r) \nonumber\\
    & \le \sum_{i=1}^{d}  \vol_{d-i}(\CC) (2\pi\g)^{\frac{i}{2}} \l(   \sqrt{2\g}\cdot \sqrt{\frac{i+3}{2}} + R  \r)
    \qquad \text{(Gautschi's inequality)}
    \nonumber\\
    & =  \sum_{i=1}^{d}  \vol_{d-i}(\CC) (2\pi\g)^{\frac{i}{2}} \l(   \sqrt{\g ( i+3)} + R  \r)\nonumber\\
    & = \sum_{i=0}^{d-1} \vol_i(\CC) (2\pi\g)^{\frac{d-i}{2}} \l( \sqrt{\g(d-i+3)}+R \r) \nonumber\\
    & =: I_2(\g),
\end{align}
which establishes \eqref{eq:int2}. 
This completes the proof of Lemma \ref{lem:integral}.

\section{Proof of Lemma \ref{lem:2ndLeg}}

Recall from \eqref{eq:Q0} that $\Q\propto e^{-f-1_\CC}$, where $1_\CC$ denotes the indicator function on the set $\CC$.  Recall also from \eqref{eq:Qgamma} that  $\Q_\g\propto e^{-f-1_{\CC,\g}}$, where $1_{\CC,\g}(x)= \frac{1}{2\g} \dist(x,\CC)^2$ is the MY envelope of $1_\CC$, see \eqref{eq:MYenvOfInd}. We will repeatedly use these two distributions in the proof. 
To begin, let us invoke Theorem~\ref{thm:WtoTV} to write that 
\begin{align}
    \W_1(\Q_\g,\Q) & \le \int \|x\|_2\cdot  |\Q_\g(x) - \Q(x) | \, \der x 
    \qquad \text{(Theorem \ref{thm:WtoTV})}
    \nonumber\\
    & = \int_{\CC} \|x\|_2 \cdot |\Q_\g(x)-\Q(x)| \, \der x + \int_{ \CC^c} \|x\|_2 \cdot  |\Q_\g (x) - \Q(x)| \, \der x.
    \label{eq:insideOutsideK}
\end{align}
For the first integral in the last line above, we use the fact that $\CC\subset \B(0,R)$ in Assumption \ref{assumption:fg}\ref{item:CAssumption} and write that
\begin{align}
    &\int_{ \CC} \|x\|_2 \cdot |\Q_\g(x)-\Q(x)| \, \der x \nonumber\\
    & \le R \int_{\CC}  |\Q_\g(x)-\Q(x)| \, \der x \qquad \text{(Assumption \ref{assumption:fg}\ref{item:CAssumption})} \nonumber\\
    & = R \int_{\CC} \l|\frac{e^{-f(x)}}{\int_{\R^d} e^{-f(z)-\frac{\dist(z,\CC)^2}{2\g}} \, \der z} - \frac{e^{-f(x)}}{\int_{\CC} e^{-f(z)}\, \der z }\r|\, \der x. 
    \nonumber\\
    & = R\int_\CC e^{-f(x)}\, \der x \l|\frac{1}{\int_{\R^d} e^{-f(z) - \frac{\dist(z,\CC)^2}{2\g}} \, \der z} - \frac{1}{\int_\CC e^{-f(z)}\,\der z}\r|
    \nonumber\\
    & = R \int_\CC e^{-f(x)}\,\der x\l(\frac{1}{\int_\CC e^{-f(z)}\,\der z} - \frac{1}{\int_{\R^d} e^{-f(z)-\frac{\dist(z,\CC)^2}{2\g}}\,\der z } \r) ,
    \label{eq:brkIntPre}
\end{align}    
where the last line uses the fact that 
$$
\int_{\R^d} e^{-f(z) - \frac{\dist(z,\CC)^2}{2\g}} \, \der z\ge 
\int_{\CC} e^{-f(z) - \frac{\dist(z,\CC)^2}{2\g}} \, \der z = \int_{\CC} e^{-f(z)}\, \der z. 
$$
We continue to simplify the last line of \eqref{eq:brkIntPre} as 
\begin{align}    
    & \int_\CC \|x\|_2\cdot |\Q_\g(x)-\Q(x)|\,\der x\nonumber\\
    & = R\l(1 - \frac{\int_\CC e^{-f(x)}\,\der x }{\int_{\R^d} e^{-f(x) - \frac{\dist(x,\CC)^2}{2\g}} \, \der x} \r) \nonumber\\
    & = R \cdot \frac{\int_{\R^d} e^{-f(x)-\frac{\dist(x,\CC)^2}{2\g}}\, \der z -  \int_\CC e^{-f(x)}\,\der x }{\int_{\R^d} e^{-f(x) - \frac{\dist(x,\CC)^2}{2\g}} \, \der x}  \nonumber\\
    & = R \cdot \frac{\int_{\CC^c} e^{-f(x)-\frac{\dist(x,\CC)^2}{2\g}}\, \der x  }{\int_{\R^d} e^{-f(x) - \frac{\dist(x,\CC)^2}{2\g}} \, \der x}   \nonumber\\
    & \le Re^{\max_x f(x) - \min_x f(x)} \frac{\int_{\CC^c} e^{-\frac{\dist(x,\CC)^2}{2\g}}\, \der x }{\int_{\R^d} e^{ - \frac{\dist(x,\CC)^2}{2\g}} \, \der x} \nonumber\\
    & = \frac{Re^{\max_x f(x) - \min_x f(x)} I_1(\g)}{\vol(\CC)+I_1(\g)}.
    \qquad \text{(Lemma \ref{lem:integral})}
        \label{eq:brkInt1}
\end{align}
For the second integral in the last line of \eqref{eq:insideOutsideK}, we can again use the definitions of $\Q$ and $\Q_\g$ to write that 
\begin{align}
    & \int_{\CC^c} \|x\|_2 \cdot |\Q_\g(x)-\Q(x)| \, \der x \nonumber\\
    & = \int_{\CC^c} \|x\|_2\cdot  \Q_\g(x) \, \der x
    \nonumber\\
    & = \frac{\int_{\CC^c} \|x\|_2 e^{-f(x)-\frac{\dist(x,\CC)^2}{2\g}} \, \der x}{\int_{\R^d} e^{-f(x)-\frac{\dist(x,\CC)^2}{2\g}} \, \der x } \nonumber\\
    & \le  \frac{e^{\max_x f(x)-\min_x f(x)} I_2(\g)}{\vol(\CC)+I_1(\g)},
    \qquad \text{(Lemma \ref{lem:integral})}
    \label{eq:brkInt2}
\end{align}
where the first identity above uses the fact that $\Q$ is supported on $\CC$. 
With \eqref{eq:brkInt1} and \eqref{eq:brkInt2} at hand, we revisit \eqref{eq:insideOutsideK} and write that
\begin{align}
     \W_1(\Q_\g,\Q) 
    & 
    \le 
    e^{\max_x f(x)-\min_x f(x)}\frac{R I_1(\g)+ I_2(\g)}{\vol(\CC)+I_1(\g)},
\end{align}
where we used \eqref{eq:insideOutsideK}, \eqref{eq:brkInt1}, and \eqref{eq:brkInt2}.
This completes the proof of Lemma~\ref{lem:2ndLeg}.

\section{Proof of Lemma \ref{lem:3rdLeg}}

Recall from  \eqref{eq:Q0} and \eqref{eq:pi}  that $\Q\propto e^{-f-1_{\CC}}$ and $\P\propto e^{-F}$, respectively.  The function $F$ was defined in~\eqref{eq:BigF}. 
In order to bound $\W_1(\Q,\P)$, our plan is  to invoke Lemma \ref{lem:h1h2lemma} with $h_1=f+1_\CC$ and $h_2=F$. As before, this means that we need to verify~\eqref{eq:h1h2} for the choice of $h_1=f+1_\CC$ and $h_2=F^\FB_\g$. We begin by relating these two functions together as 
\begin{align}
    F(x) & = f(x)+g(x) =f(x)+g(x)+1_{\CC}(x) \qquad \text{(see \eqref{eq:BigF})} \nonumber\\
    & \le \max_{z\in \CC} g(z) +f(x) + 1_\CC(x).
\end{align}
In the first line above, we used the fact that  $F$ and $1_\CC$ both take infinity outside of the set  $\CC$. 
In the other direction, we write that 
\begin{align}
    F(x) & = f(x)+g(x)+1_\CC(x) \nonumber\\
    & \ge \min_{z\in \CC}g(z)+ f(x) + 1_\CC(x).
\end{align}
To summarize, for our choice of $h_1=f+1_\CC$ and $h_2=F$, \eqref{eq:h1h2} is satisfied with
\begin{align}
    & \b' (f(x)+1_\CC(x))+\b \le F(x)\le (f(x)+1_\CC(x))+\a, \qquad x\in \R^d, \label{eq:summary3}\\
    & \a = \max_{z\in \CC}g(z),\quad \b = 
    \min_{z\in \CC}g(z),\quad \b' = 1.\nonumber
\end{align}
We can now invoke Lemma \ref{lem:h1h2lemma} to find that 
\begin{align}
    \W_1(\Q,\P ) & = \W_1\l(\frac{e^{-f-1_{\CC}}}{\int e^{-f-1_\CC}} ,
    \frac{e^{-F}}{\int e^{-F}}
    \r) \nonumber\\
    & \le  \l(e^{\max_{x\in \CC}g(x)-\min_{x\in \CC}g(x) }-e^{\min_{x\in \CC}g(x)-\max_{x\in \CC} g(x)}\r) \nonumber\\
    & \qquad \cdot 
    \frac{\int \|x\|_2 e^{-f(x)-1_\CC(x)} \, \der x}{
    \int e^{-f(x)-1_\CC(x)} \, \der x
    } \qquad \text{(Lemma \ref{lem:h1h2lemma})} \nonumber\\
    & \le  \l(e^{\max_{x\in \CC}g(x)-\min_{x\in \CC}g(x) }-e^{\min_{x\in \CC}g(x)-\max_{x\in \CC} g(x)}\r) \nonumber\\
    &\qquad     \cdot \frac{R \int e^{-f(x)-1_{\CC}(x)}\,\der x}{\int e^{-f(x)-1_{\CC}(x)}\,\der x}
    \qquad (\CC \subset \B(0,R)) \nonumber\\
    & =  \l(e^{\max_{x\in \CC}g(x)-\min_{x\in \CC}g(x) }-e^{\min_{x\in \CC}g(x)-\max_{x\in \CC} g(x)}\r) R,
\end{align}
where we used the fact that $\CC \subset \B(0,R)$ by Assumption \ref{assumption:fg}\ref{item:CAssumption} in the penultimate line above. This completes the proof of Lemma \ref{lem:3rdLeg}.

\section{Proof of Lemma \ref{lem:termTwo}}

Consider the stochastic differential equation
\begin{equation}
    \der x_t = -\nabla F_{\g_{k}}(x_t)+\sqrt{2}\der B_t,
    \qquad x_0 \sim \P_{\g_{k}},
    \label{eq:sde}
\end{equation}
where $\{B_t\}_{t\ge 0}$ is the standard Brownian motion. Above, note that the initial probability measure is~$\P_{\g_k}\propto e^{-F_{\g_k}}$.  From Definition~\ref{assump:envAssumption}\ref{item:gradLips}, recall that $F_{\g_{k}}$ is $\lambda_{\g_{k}}$-smooth and is, moreover, coercive. To be concrete, the latter means that $F_{\g_{k}}$
satisfies $$\lim_{\|x\|_2\rightarrow\infty} F_{\g_{k}}(x)=\infty.$$ Therefore Theorem~3.4 in \cite{eberle} ensures that $\P_{\g_{k}}\propto e^{-F_{\g_{k}}}$ is the invariant probability measure of \eqref{eq:sde}. In particular, because $x_0\sim \P_{\g_{k}}$ in \eqref{eq:sde}, it holds that
\begin{equation}
    x_t \sim \P_{\g_{k}}, \qquad t\ge 0.
    \label{eq:remainsInv}
\end{equation}
Let $\{P_{k}\}_{k\ge 0}$ denote the Markov transition kernel associated with the Markov chain $\{x_{H_k}\}_{k\ge 0}$, where 
\begin{equation}
H_k := \sum_{i=1}^{k-1} h_i
\label{eq:Hk}
\end{equation}
is the elapsed time since initialization. Using this transition kernel, we can rewrite \eqref{eq:remainsInv} as 
\begin{equation}
    \P_{\g_{k}} =  \P_{\g_{k}} P_{k}, \qquad k\ge 0.
    \label{eq:Pp}
\end{equation}
Finally, we can write the quantity of interest on the left-hand side of \eqref{eq:discErr} as 
\begin{align}
    \W_1( \P_{\g_{k}} Q_{k} , \P_{\g_{k}}) =  \W_1( \P_{\g_{k}} Q_{k}, \P_{\g_{k}} P_{k}). 
    \qquad \text{(see \eqref{eq:Pp})}
    \label{eq:W1Rewritten}
\end{align}
In the remainder of the proof, we will upper bound the right-hand side above, which can be thought of as the discretization error associated with~\eqref{eq:sde}. To begin, recall from \eqref{eq:remainsInv} that $x_{H_{k}}\sim \P_{\g_{k}}$, and note that 
\begin{align}
x'& :=x_{H_k}-h_k \nabla F_{\g_k}(x_{H_k})+\sqrt{2h_k}\zeta_{k+1} \sim \P_{\g_{k}} Q_{k} 
\qquad \text{(see \eqref{eq:QTransitionKernel})}
\nonumber\\
& = x_{H_k} - \int_{H_k}^{H_{k+1}} \nabla F_{\g_k}(x_{H_k}) \der t + \sqrt{2h_k}\zeta_{k+1},\qquad \text{(see \eqref{eq:Hk})}
\label{eq:1stLegofW1}
\end{align} 
by construction. The second line above uses the fact that $H_{k+1}-H_k=h_k$ by~\eqref{eq:Hk}. 
Likewise, note that 
\begin{align}
& x_{H_{k+1}} \nonumber\\
& =x_{H_k}-\int_{H_k}^{H_{k+1}} \nabla F_{\g_{k}}(x_t) \der t + \sqrt{2}\int_{H_k}^{H_{k+1}} \der B_t \sim  \P_{\g_{k}} P_{k}
\nonumber\\
& \overset{\text{dist.}}{=} x_{H_k} -\int_{H_k}^{H_{k+1}} \nabla F_{\g_{k}}(x_t) \der t + \sqrt{2h_{k}}\zeta_{k+1}, 
\label{eq:2ndLegofW1}
\end{align}
where the first line uses  \eqref{eq:sde} and \eqref{eq:Pp}, and 
both sides of $\overset{\text{dist.}}{=}$ have the same distribution. Above, recall that $\zeta_{k+1}\sim \text{normal}(0,I_d)$.  
In view of~\eqref{eq:1stLegofW1} and~\eqref{eq:2ndLegofW1}, we revisit \eqref{eq:W1Rewritten} and write that 
\begin{align}
    & \W_1(\P_{\g_{k}} Q_{k} , \P_{\g_{k}}) \nonumber\\
    & =  \W_1( \P_{\g_{k}} Q_{k}, \P_{\g_{k}} P_{k})
    \qquad \text{(see \ref{eq:W1Rewritten})}\nonumber\\
    & \le \E \| x'-x_{H_{k+1}}\|_2 
    \qquad \text{(see \eqref{eq:wassDist}, \eqref{eq:1stLegofW1},  \eqref{eq:2ndLegofW1})}
    \nonumber\\
    & = \E \l\| \int_{H_k}^{H_{k+1}} \nabla F_{\g_{k}}(x_t) - \nabla F_{\g_{k}}(x_{H_k}) \der t \r\|_2
    \qquad \text{(see \eqref{eq:1stLegofW1},  \eqref{eq:2ndLegofW1})}
    \nonumber\\
    & \le \int_{H_k}^{H_{k+1}} \E \|F_{\g_{k}}(x_t) - \nabla F_{\g_{k}}(x_{H_k}) \|_2 \der 
    t
    \qquad \text{(Minkowski's integral inequality)}
    \nonumber\\
    & \le \lambda_{\g_{k}} \int_{H_k}^{H_{k+1}} \E \| x_t - x_{H_k}\|_2 \der t ,
    \qquad \text{(Definition \ref{assump:envAssumption}\ref{item:gradLips})}
    \label{eq:W1toE}
\end{align}
where the last line above uses the fact that $F_{\g_{k}}$ is $\lambda_{\g_{k}}$-smooth.
To estimate the expectation in the last line above, we write that 
\begin{align}
    & \E \|x_t - x_{H_k}\|_2 \nonumber\\
    & = \E \l\| \int_{H_k}^t \der x_s \r\|_2 \nonumber\\
    & =  
    \E \l\| - \int_{H_k}^t  \nabla F_{\g_{k}}(x_s) \der s + \sqrt{2} \int_{H_k}^t \der B_s  \r\|_2
    \qquad \text{(see \eqref{eq:sde})} \nonumber\\
    & \le \E \l\| \int_{H_k}^t \nabla F_{\g_{k}} (x_s) \der s \r\|_2 + \sqrt{2} \E \l\| \int_{H_k}^t \der B_s\r\|_2 
    \qquad \text{(triangle inequality)}
    \nonumber\\
    & \le \int_{H_k}^t \E \| \nabla F_{\g_{k}} (x_s) \|_2 \der s +\sqrt{2} \E \l\| \int_{H_k}^t \der B_s\r\|_2.
    \, \text{(Minkowski's inequality)}
    \label{eq:brkTwo}
\end{align}
To estimate the first expectation in the last line of \eqref{eq:brkTwo}, recall from the last line of \eqref{eq:W1toE} that $t\in [H_k,H_{k+1})$ and then note that 
\begin{align}
& \int_{H_k}^t \E \| \nabla F_{\g_{k}} (x_s) \| \der s \nonumber\\ & =  \int_{H_k}^t \E_{x\sim \P_{\g_{k}}} \| \nabla F_{\g_{k}} (x) \|_2 \der s   \qquad \text{(see \eqref{eq:remainsInv})}\nonumber\\
    & = (t-H_k) \E_{x\sim \P_{\g_{k}}} \| \nabla F_{\g_{k}} (x) \|_2 \nonumber\\
    & \le h_k \E_{x\sim \P_{\g_{k}}} \| \nabla F_{\g_{k}} (x) \|_2 
    \qquad \text{(see \eqref{eq:Hk})}
    \nonumber\\
    & \le h_k  \sqrt{\E_{x\sim \P_{\g_{k}}} \| \nabla F_{\g_{k}} (x)\|_2^2  } 
    \qquad \text{(Jensen's inequality)}
    \nonumber\\
    & = h_k \sqrt{\E_{x\sim \P_{\g_{k}}} \l[\Delta F_{\g_{k}}(x)\r] }
    \qquad (\text{integration by part and } \P_{\g_{k}}\sim e^{-F_{\g_{k}}})
    \nonumber\\
    & \le h_k \sqrt{  \lambda_{\g_{k}} d}, 
    \label{eq:brk1}
\end{align}
where the last line above follows from $\Delta F_{\g_{k}}(x)=\tr(\nabla^2 F_{\g_{k}}(x))$ and the fact that $F_{\g_{k}}$ is $\lambda_{\g_k}$-smooth by Definition~\ref{assump:envAssumption}\ref{item:gradLips}.
To estimate the second expectation in the last line of \eqref{eq:brkTwo}, note that 
\begin{equation}
\int_{H_k}^t \der B_s \overset{\text{dist.}}{=} \sqrt{t-H_k} \zeta_{k+1}, \qquad \zeta_{k+1}\sim \text{normal}(0,I_d).
\label{eq:propsBrownian}
\end{equation}
The above observation enables us to bound the second expectation in the last line of \eqref{eq:brkTwo} as 
\begin{align}
\E \l\| \int_{H_k}^t \der B_s\r\|_2& 
    = \sqrt{t-H_k} \E \| \zeta_{k+1}\|_2 
    \qquad \text{(see \eqref{eq:propsBrownian})}
    \nonumber\\
    & \le \sqrt{t-H_k} \sqrt{\E\| \zeta_{k+1}\|_2^2} 
    \qquad \text{(Jensen's inequality)}
    \nonumber\\
    & = \sqrt{(t-H_k)d}  \qquad \text{(see \eqref{eq:propsBrownian})} \nonumber\\
    & \le \sqrt{ h_{k} d}. 
    \qquad \text{(see \eqref{eq:Hk} and \eqref{eq:W1toE})}
    \label{eq:brk2}
\end{align}
We now plug in the bounds in \eqref{eq:brk1} and \eqref{eq:brk2} back into \eqref{eq:brkTwo} to obtain that 
\begin{align}
  \E \|x_t - x_{H_k}\|_2 
    & \le \int_{H_k}^t \E \| \nabla F_{\g_{k}} (x_s) \|_2 \der s +\sqrt{2} \E \l\| \int_{H_k}^t \der B_s\r\|_2 \qquad \text{(see \eqref{eq:brkTwo})} \nonumber\\
    & \le \sqrt{h_k d} \l( \sqrt{h_k \lambda_{\g_{k}} } +\sqrt{2}\r). 
    \qquad \text{(see \eqref{eq:brk1} and \eqref{eq:brk2})}
    \label{eq:fixedT}
\end{align}
By substituting the above bound back into \eqref{eq:W1toE}, we arrive at
\begin{align}
    \W_1( \P_{\g_{k}} Q_{k}, \P_{\g_{k}})& 
    \le \lambda_{\g_{k}} \int_{H_k}^{H_{k+1}} \E \| x_t - x_{H_k}\|_2 \der t 
    \qquad \text{(see \eqref{eq:W1toE})}
    \nonumber\\
    & \le  \lambda_{\g_{k}}(H_{k+1}-H_k) \sqrt{h_k d} \cdot \l( \sqrt{h_k \lambda_{\g_{k}} } +\sqrt{2}\r)
    \qquad \text{(see \eqref{eq:fixedT})} \nonumber\\
    & = \lambda_{\g_{k}}h_k \sqrt{h_k d} \cdot \l( \sqrt{h_k \lambda_{\g_{k}} } +\sqrt{2}\r) \qquad \text{(see \eqref{eq:Hk})} \nonumber\\
    & = \lambda_{\g_{k}} \sqrt{h_k^3 d} \cdot \l( \sqrt{h_k \lambda_{\g_{k}} } +\sqrt{2}\r),
\end{align}
which completes the proof of Lemma \ref{lem:termTwo}.

\end{document}